\documentclass{article}   

\usepackage{graphicx}
\usepackage{epsfig,mathrsfs,amsmath,amsfonts,amssymb}
\usepackage{color}
\usepackage{multirow}
\usepackage{tabularx}

\usepackage{geometry}
\usepackage{authblk}

\usepackage{pgf}
\usepackage{tikz}

\providecommand{\keywords}[1]{\textbf{\textit{Keywords:}} #1}
\newtheorem{remark}{Remark}
\newtheorem{example}{Example}
\newcommand\preprint[1]{
  \begingroup
  \renewcommand\thefootnote{}\footnote{\\#1}
  \addtocounter{footnote}{-1}
  \endgroup
}

\begin{document}

\title{Improved Accuracy of High-Order WENO Finite Volume Methods on
  Cartesian Grids\thanks{This work was supported by the DFG through FOR1048.}
}

\author{Pawel Buchm\"uller\thanks{Pawel.Buchmueller@rub.de}}
\author{Christiane Helzel\thanks{Corresponding author. Christiane.Helzel@rub.de}}
\affil{Department of Mathematics, Ruhr-University Bochum}

\date{January 29, 2014}

\maketitle 

\preprint{Preprint submitted to {\em  Journal of Scientific Computing}. \\
  The final publication is available at 
  \href{http://dx.doi.org/10.1007/s10915-014-9825-1}{\em link.springer.com}}

\begin{abstract}
We propose a simple modification of standard WENO finite volume methods for
Cartesian grids, which retains the full spatial order of accuracy 
of the one--dimensional discretization when applied to
nonlinear multidimensional systems of conservation laws.

We derive formulas, which allow us to compute high-order accurate
point values of the conserved quantities at grid cell
interfaces. Using those point values, we can compute a high-order flux
at the center of a grid cell interface. Finally, we use those point
values to compute high-order accurate averaged fluxes at cell
interfaces as needed by a finite volume method. 

The method is described in detail for the two--dimensional Euler equations of gas
dynamics. An extension to the three--dimensional case as well as
to other nonlinear systems of conservation laws in divergence form is straightforward. 
Furthermore, similar ideas can be used to improve the accuracy of WENO
type methods for hyperbolic systems which are not in divergence form.

Several test computations 
confirm the high-order accuracy for smooth nonlinear problems. \\

\end{abstract}

\section{Introduction}
High--order WENO (i.e., weighted essentially non-oscillatory) methods are wi\-dely used  for the approximation of
hyperbolic problems, see for example the recent review of Shu \cite{article:Shu2009}. 
The simplest way to use WENO methods
on multidimensional Cartesian grids consists in applying a
one--dimensional WENO scheme in each direction. This spatial
discretization is typically
combined with a Runge--Kutta method in time, i.e.\ during each stage of
the Runge--Kutta method one--dimensional WENO schemes are used in a
dimension--by--dimension fashion.

On uniform Cartesian grids, conservative finite difference WENO
methods based on flux interpolation, as introduced by Shu and Osher
\cite{article:SO88,article:SO89}, lead to high order accurate
approximations of the conserved quantities for linear as well as
nonlinear conservation laws. An extension to smoothly varying mapped
grids is possible, see \cite{article:Merriman2003}.
In contrast to this, finite volume WENO methods based on a 
dimension--by--dimension approach   retain the full order of accuracy for smooth solutions of
linear multi--dimensional problems but they are only second order
accurate for smooth solutions of nonlinear problems, see
\cite{article:Shu2009,article:ZZS2011}.

Here we restrict our considerations to finite volume WENO methods.
For hyperbolic equations in divergence form, an advantage of finite volume methods
is that they approximate the integral form of a conservation law which remains
valid at discontinuities, where the differential form of the equation
is not valid in the classical sense. In general it is straight forward
to extend finite volume methods to unstructured grids. With a 
dimension--by--dimension approach we are of course limited to
Cartesian grids, but not necessarily to equidistant Cartesian grids. 
The standard approach to avoid the loss of accuracy of the
dimension--by--dimension approach is to use a multidimensional
reconstruction (with WENO limiting) and a high order quadrature
formula to compute fluxes at grid cell interfaces.
Such methods were used on unstructured
as well as on structured grids in \cite{article:CA1993,article:TTD2011,article:ZZS2011,article:HRT2013}. 
However, the multidimensional reconstruction and in particular the
limiting of the reconstructed polynomials is quite expensive.
Furthermore, such schemes require  flux computations at several points
per grid cell interface (i.e., at the nodes of the quadrature
formula). 

A method first described in the context of ENO methods in
\cite{article:CA1993} and later used for the construction of WENO
methods on Cartesian grids in \cite{article:SHS2002,article:TT2004,article:ZZS2011} is most closely related to our
approach in the sense that it aims to overcome the formal loss of
accuracy without using a full multi--dimensional polynomial
reconstruction of the conserved quantities. Instead, those authors use
one--dimensional WENO reconstruction
to obtain a high--order accurate approximation of face--averaged
values of the conserved quantities at all grid cell interfaces. In a
second reconstruction step, these averaged values at grid cell interfaces 
are used to construct a one--dimensional polynomial representation
of the conserved quantities along grid cell interfaces. 
Thus, for each grid cell interface, a one dimensional WENO reconstruction is
required in the $x$ as well as the $y$--direction.  
Finally, the method of \cite{article:ZZS2011} computes fluxes at grid
cell interfaces by using a high--order accurate quadrature formula,
thus it requires the evaluation of the numerical flux at all nodes of
the quadrature formula.

Here we present a simpler modification of finite volume WENO methods, 
which also leads to the full spatial order of
accuracy by using only one--dimensional polynomial
reconstructions in a dimension--by--dimension approach.
While WENO reconstruction is typically of odd order (here we consider
methods of order five and seven), the corrections introduced in this
paper lead to fluxes of even order (here we present the formulas for
order  four and six).
For the temporal discretization we use explicit Runge-Kutta methods of
order five or seven.
An important component of our approach is the transfer of high--order
averaged values of the conserved quantities to high--order point values
and vice versa. 
Such a transformation has also been used  
in the recently proposed fourth order accurate finite volume method of
McCorquodale and Colella \cite{article:MC2011}. 

In this paper we restrict our considerations to finite volume WENO
methods on equidistant Cartesian grids, i.e.\ a case that can be
handled perfectly well by finite difference WENO methods. Even on such
grids, there are situations where finite volume methods are more
appropriate than finite difference methods. For example, if we want to
construct a high order accurate and conservative method that uses
adaptive mesh refinement (AMR), see
\cite{article:SQC2011} for more discussions. 

\section{The dimension--by--dimension 
WENO finite volume method for hyperbolic problems in divergence form}
\label{section:s2}
In this section we give a brief description of Cartesian grid WENO methods
and review their accuracy for linear and nonlinear multidimensional
problems.

We consider two--dimensional systems of conservation laws, i.e.\
initial value problems of the form 
\begin{equation}\label{eqn:2dproblem}
\begin{split}
\partial_t q + \partial_x f(q) + \partial_y g(q) = 0\\
q(x,y,0) = q_0(x,y)
\end{split}
\end{equation}
where $q:\mathbb{R}^2 \times \mathbb{R}^+ \rightarrow \mathbb{R}^m$ is
a vector of conserved quantities, and $f,\ g : \mathbb{R}^m \rightarrow
\mathbb{R}^m$ are vector valued flux functions. 

In order to discretize (\ref{eqn:2dproblem}), we use a method of lines
approach. 
We restrict our considerations to equidistant Cartesian grids with
grid cells $C_{i,j} = (x_{i- \frac{1}{2}},x_{i+\frac{1}{2}}) \times
(y_{j-\frac{1}{2}},y_{j+\frac{1}{2}})$ and mesh width $\Delta x =
x_{i+\frac{1}{2}}-x_{i-\frac{1}{2}}$ and $\Delta y =
y_{j+\frac{1}{2}}-y_{j-\frac{1}{2}}$ for all $i,j$.
A  finite volume method can be written in the semi--discrete form
\begin{equation}\label{eqn:semi-discrete-fv}
\frac{d}{dt} Q_{i,j}(t) = - \frac{1}{\Delta x} \left(
  F_{i+\frac{1}{2},j}(t) - F_{i-\frac{1}{2},j}(t) \right) -
\frac{1}{\Delta y} \left( G_{i,j+\frac{1}{2}} (t) -
  G_{i,j-\frac{1}{2}}(t) \right),
\end{equation}
where $Q_{i,j}(t)$ is an approximation of the cell average of the conserved quantities in
grid cell $C_{i,j}$ and the terms $F(t),\ G(t)$ are flux functions at the
grid cell interfaces in the $x$ and the $y$--direction, respectively.
For the temporal discretization we use explicit Runge--Kutta methods 
of appropriate order of accuracy. The two different Runge-Kutta
methods, RK5 and RK7, used for
our computations are described in Appendix \ref{app:RK4-6}.

For the spatial discretization we use in each direction a one--dimensional piecewise polynomial
reconstruction of the conserved quantities. In the $x$--direction
we construct one--dimensional polynomials $q_{i,j}^1(x)$ and in the
$y$--direction we construct polynomials $q_{i,j}^2(y)$. These
polynomials are local approximations of the conserved quantity in cell
$C_{i,j}$. Furthermore, they satisfy
\begin{equation}
\begin{split}
Q_{i,j} & = 
\frac{1}{\Delta x} \int_{x_{i-\frac{1}{2}}}^{x_{i+\frac{1}{2}}}
q_{i,j}^1(x) dx\\
& = \frac{1}{\Delta y} \int_{y_{j-\frac{1}{2}}}^{y_{j+\frac{1}{2}}}
q_{i,j}^2(y) dy.
\end{split}
\end{equation}

For each grid cell of a two--dimensional Cartesian grid, we reconstruct
four edge averaged values of the
conserved quantities by evaluating these
polynomials at the interfaces. Those interface values of the conserved
quantities are denoted by
\begin{equation}
\begin{split}
& Q_{i-\frac{1}{2},j}^+ := 
q_{i,j}^1(x_{i-\frac{1}{2}}), \quad Q_{i+\frac{1}{2},j}^- := q_{i,j}^1(x_{i+\frac{1}{2}}), \\
& Q_{i,j-\frac{1}{2}}^+ := 
q_{i,j}^2(y_{j-\frac{1}{2}}), \quad Q_{i,j+\frac{1}{2}}^- :=  q_{i,j}^2(y_{j+\frac{1}{2}}).
\end{split}
\end{equation} 
Here we use component--wise WENO reconstruction
of order five and seven, known as WENO-Z method, see Appendix \ref{appendix:s2}.
Note that at each grid cell interface we have two
reconstructed edge--averaged values of the conserved quantities.
Assuming that the WENO reconstruction was based on exact cell average
values of the conserved quantities, then the edge averaged values
satisfy (for sufficiently smooth functions $q$) 
\begin{equation}\label{eqn:interface}
  \begin{aligned}
    Q_{i-\frac{1}{2},j}^\pm (t) & = \frac{1}{\Delta y}
    \int_{y_{j-\frac{1}{2}}}^{y_{j+\frac{1}{2}}} q(x_{i-\frac{1}{2}},y,t) dy
    +\mathcal{O}(\Delta x^p)\\
    Q_{i,j-\frac{1}{2}}^\pm (t) & = \frac{1}{\Delta x} 
    \int_{x_{i-\frac{1}{2}}}^{x_{i+\frac{1}{2}}} q(x,y_{j-\frac{1}{2}},t)
    dx 
    +\mathcal{O}(\Delta y^p)
  \end{aligned}
\end{equation}
with $p=5$ or $p=7$, respectively.
In general, the cell averaged values are only $p$-th order accurate
approximations of the exact cell averaged values, i.e.\
\begin{equation}
Q_{i,j}(t) = \frac{1}{\Delta x \Delta y}
\int_{y_{j-\frac{1}{2}}}^{y_{j+\frac{1}{2}}}
\int_{x_{i-\frac{1}{2}}}^{x_{i+\frac{1}{2}}} q(x,y,t) dx dy +\mathcal{O}
(\Delta x^p + \Delta y^p). 
\end{equation}
We then get, for both equations of   (\ref{eqn:interface}), on the
right hand side an  error of the form ${\cal O}(\Delta x^p + \Delta y^p)$.

The numerical fluxes 
$F_{i-\frac{1}{2},j}(t)$ 
and
$G_{i,j-\frac{1}{2}}(t)$
can be obtained by using a numerical flux function such as 
Lax-Friedrichs, which has the form
\begin{equation}\label{eqn:LxF}
F_{i-\frac{1}{2},j} = \frac{1}{2} \left[ f(Q_{i-\frac{1}{2},j}^- ) +
  f(Q_{i-\frac{1}{2},j}^+) - \alpha (Q_{i-\frac{1}{2},j}^+ -
  Q_{i-\frac{1}{2},j}^-) \right],
\end{equation}
where $\alpha$ is an upper estimate for the largest absolute value of the
eigenvalues of the flux Jacobian matrix. 

Alternatively, we can compute an interface value $Q_{i-\frac{1}{2},j}^*$
of the conserved quantities, by exact or approximative solution of the
Riemann problem 
with data $Q_{i-\frac{1}{2},j}^\pm$. The flux  can then be computed using $F_{i-\frac{1}{2},j} =
f(Q_{i-\frac{1}{2},j}^*)$. 

\begin{remark}
For the computation of smooth solution of the Euler equation, the
choice of the numerical flux function has only a very small effect on
the quality of the numerical solutions in high-order WENO methods.
Therefore, we used the
Lax-Friedrichs flux function for the convergence studies shown in
Section \ref{section:results}.
For the computation of problems with shock waves or contact
discontinuities, the choice of the
numerical flux function has a larger impact on the quality of the
numerical solution. For such problems we used the Roe Riemann solver,
with an entropy fix according to Harten and Hyman,
in order to compute the conserved quantity at the grid cell interfaces
and evaluate the flux for this value.
\end{remark}

In the linear case, i.e.\  for $f(q) = A q$ with a constant matrix $A
\in {\mathbb R}^{m\times m}$, we obtain (by using the Lax--Friedrichs flux)
\begin{equation*}
\begin{split}
F_{i-\frac{1}{2},j} & = \frac{1}{2} \left[ A(Q_{i-\frac{1}{2},j}^- +
  Q_{i-\frac{1}{2},j}^+) - \alpha (Q_{i-\frac{1}{2},j}^+ -
  Q_{i-\frac{1}{2},j}^-) \right] \\
& = \frac{1}{2}  A  \frac{2}{\Delta y}
\int_{y_{j-\frac{1}{2}}}^{y_{j+\frac{1}{2}}} q(x_{i-\frac{1}{2}},y,t)
dy +\mathcal{O}(\Delta x^p + \Delta y^p) \\
& = \frac{1}{\Delta y} \int_{y_{j-\frac{1}{2}}}^{y_{j+\frac{1}{2}}} A
q(x_{i-\frac{1}{2}},y,t) dy +\mathcal{O}(\Delta x^p + \Delta y^p),
\end{split}
\end{equation*}
and thus a $p$-th order accurate approximation of the average value of
the flux across the interface.
This is exactly what is needed in order to construct a $p$-th order
accurate finite volume scheme.
Interfaces in the $y$--direction are treated
analogously. 

In the nonlinear case, the dimension--by--dimension approach is in
general only second order accurate. The reason for this loss of
accuracy is that
a flux function applied to edge averaged values of the conserved
quantities does not provide an averaged flux of the same accuracy.  
This explains the different accuracy of multi--dimensional finite
difference and finite volume WENO methods, see
\cite{article:Shu2009,article:ZZS2011}.
While finite difference WENO methods retain the full order of
accuracy of the one--dimensional reconstruction, finite volume WENO
methods are in general only second order accurate.

The $p$-th order accurate edge average value $Q_{i-\frac{1}{2},j}^\pm(t)$
provides a second order accurate approximation of the point value of
the conserved quantities at the midpoint of the cell edge, i.e.\
\begin{equation}\label{eqn:midpoint}
\begin{split}
Q_{i-\frac{1}{2},j}^\pm(t) & = \frac{1}{\Delta y}
\int_{y_{j-\frac{1}{2}}}^{y_{j+\frac{1}{2}}} q(x_{i-\frac{1}{2}}, y,
t) dy +\mathcal{O}(\Delta x^p)\\
&= q(x_{i-\frac{1}{2}},y_{j},t) + {\cal
  O}(\Delta y^2 + \Delta x^p).
\end{split}
\end{equation} 
This can be interpreted as using the midpoint rule to approximate the
integral on the right hand side of (\ref{eqn:midpoint}).
Using the Lax-Friedrichs flux (\ref{eqn:LxF}) and assuming
that the flux function $f$ can be expanded using Taylor series expansion, 
we get
\begin{equation}\label{eqn:aux1}
F_{i-\frac{1}{2},j}(t)  =  f(q(x_{i-\frac{1}{2}},y_j,t)) +\mathcal{O}(\Delta x^p + \Delta y^2). 
\end{equation}
Furthermore, we have
\begin{equation}\label{eqn:aux2}
f(q(x_{i-\frac{1}{2}},y_j,t)) = \frac{1}{\Delta y}
\int_{y_{j-\frac{1}{2}}}^{y_{j+\frac{1}{2}}}
f(q(x_{i-\frac{1}{2}},y,t)) dy +\mathcal{O}(\Delta y^2).
\end{equation}
From (\ref{eqn:aux1}) and (\ref{eqn:aux2}) we conclude that
\begin{equation}
F_{i-\frac{1}{2},j}(t) = \frac{1}{\Delta y}
\int_{y_{j-\frac{1}{2}}}^{y_{j+\frac{1}{2}}}
f(q(x_{i-\frac{1}{2}},y,t)) dy +\mathcal{O}(\Delta x^p + \Delta y^2).
\end{equation} 
Thus the resulting WENO method is in general only second order accurate in space.
We now summarize the
dimension--by--dimension WENO method.

{\bf Algorithm: Dimension--by--dimension WENO method.}
\nopagebreak
\begin{enumerate}
\item Compute averaged values of the conserved quantities at grid cell
  interfaces using one--dimensional WENO reconstruction, i.e.\ compute
$$
Q_{i-\frac{1}{2},j}^\pm(t), Q_{i,j-\frac{1}{2}}^\pm(t) 
$$
at all grid cell interfaces.
\item Compute fluxes at grid cell interfaces using a consistent
  numerical flux function (such as Lax-Friedrichs), i.e.\
$$
F_{i-\frac{1}{2},j}(t) = {\cal F}(Q_{i-\frac{1}{2},j}^-,
Q_{i-\frac{1}{2},j}^+), \quad
G_{i,j-\frac{1}{2}}(t) = {\cal F}(Q_{i,j-\frac{1}{2}}^-, Q_{i,j-\frac{1}{2}}^+)
$$
\item Approximate the semi--discrete system (\ref{eqn:semi-discrete-fv}), using a
  high--order accurate Runge--Kutta method.
\end{enumerate}

\section{A modification of the dimension--by--dimension approach}
We now describe a simple modification of the WENO method, which increases the accuracy
of the dimension--by--dimension approach. With this modification, full order of
accuracy can be retained for multidimensional nonlinear problems. 
The method is computationally less expensive than the method used in 
\cite{article:SHS2002,article:TT2004,article:ZZS2011}, since it is based on the one--dimensional
reconstructions used in 
the dimension--by--dimension approach.  Furthermore, our approach 
requires only one evaluation of the flux function per interface.

WENO reconstruction provides us with high order accurate
approximations of averaged values of the conserved quantities at grid
cell interfaces.
For conservation laws with nonlinear flux functions, we can not
directly compute high 
order accurate averaged values of the interface flux from these
edge averaged values of the conserved quantities.

In order to compute accurate flux functions, we first compute 
point values of the conserved quantities at the
midpoint of the grid cell interface. We then compute the numerical
flux at the point value and finally compute averaged values of the
flux at grid cell interfaces.

\subsection{Transformation between average values and point values}
We discuss the transformation between average values and point values
for functions of one spatial variable. This agrees with the situation
which will later be used in our method, since the second variable at
grid cell interfaces will just lead to an additional
index. Furthermore,  we often suppress the time dependence of the
functions in this subsection.  

We denote with $Q_i$ an approximation of the cell average of the
function $q$ in
grid cell $i$, i.e.\ the interval
$(x_{i-\frac{1}{2}},x_{i+\frac{1}{2}})$ and by $q_i$ 
an approximation of the point value $q(x_i)$ of the quantities  $q$ 
at the midpoint $x_i$ of the grid cell.
For sufficiently smooth functions $q:\mathbb{R}\rightarrow \mathbb{R}^m$,  Taylor series expansion
provides 
\begin{equation*}
\begin{split}
Q_i & = \frac{1}{\Delta x} \int_{x_{i-\frac{1}{2}}}^{x_{i+\frac{1}{2}}}
q(x) dx = \frac{1}{\Delta x} \int_{-\frac{\Delta x}{2}}^{\frac{\Delta
    x}{2}} q(x_i + x) dx \\
& = \frac{1}{\Delta x} \int_{-\frac{\Delta x}{2}}^{\frac{\Delta x}{2}}
\left( q(x_i) + x q'(x_i) + \frac{x^2}{2} q''(x_i) + \frac{x^3}{6}
  q'''(x_i) + \frac{x^{4}}{24} q^{(4)}(x_i) + \ldots \right) dx
\end{split}
\end{equation*}
and thus  the transformation 
\begin{equation}\label{eqn:transform}
 q_i = Q_i  - \frac{\Delta x^2}{24} q''(x_i) - \frac{\Delta x^4}{1920}
q^{(4)} (x_i) + \ldots 
\end{equation}
between point values and cell average values.

Thus we need expressions for the approximation of second and fourth
derivatives. In order to transform from point values to cell average values, we can
approximate these derivatives using standard finite difference
formulas.
If we transform from cell average values to point values, we use cell
average values of the conserved quantities in
order to approximate the second and fourth derivative at the midpoint
of the interval.

\subsubsection{Approximation of  derivatives from point values}
The second derivative $q^{''}(x_i)$ can be approximated using
point values of the quantity $q$ via the well known second order
accurate FD formula
\begin{equation}\label{eqn:q2_xi}
q^{''}(x_i) = \frac{1}{\Delta x^2}\left( q_{i-1} - 2 q_i +
  q_{i+1} \right) +\mathcal{O}(\Delta x^2).
\end{equation}
A fourth order accurate representation of $q^{''}(x_i)$ can be
obtained using the formula
\begin{equation}
q^{''}(x_i) = \frac{1}{12 \Delta x^2} \left( -q_{i-2} + 16 q_{i-1} -
30 q_i + 16 q_{i+1} -  q_{i+2} \right) + {\cal
  O}(\Delta x^4).
\end{equation}

A second order accurate representation of $q^{(4)}(x_i)$ can be
computed from point values using the finite difference formula
\begin{equation}\label{eqn:q4_xi}
q^{(4)}(x_i) = \frac{1}{\Delta x^4} \left(q_{i-2} - 4 q_{i-1} + 6 q_i - 4 q_{i+1} +
q_{i+2}\right) +\mathcal{O}(\Delta x^2).
\end{equation}
All of these formulas can be verified using Taylor series expansion.

\subsubsection{Approximation of derivatives from cell average values} 
Second order accurate approximations of $q''(x_i)$ and
$q^{(4)}(x_i)$ can be obtained using formulas analogously to (\ref{eqn:q2_xi}) and
(\ref{eqn:q4_xi}) with the point values $q_i$ replaced by average
values $Q_i$.

A fourth order accurate approximation of $q^{''}(x_i)$ can be computed
from cell averaged values via the formula
\begin{equation}\label{eqn:q2-ave}
q^{''}(x_i) = \frac{1}{8 \Delta x^2} \left(- Q_{i-2} + 12 Q_{i-1} -
22 Q_i + 12 Q_{i+1} -  Q_{i+2} \right) + {\cal
  O}(\Delta x^4).
\end{equation}
This formula can also be verified by Taylor series expansion,
after using (\ref{eqn:transform}) to
express the cell average values by point values.

Note that the equations (\ref{eqn:q2_xi})-(\ref{eqn:q2-ave}) are only valid for
uniform Cartesian grids. Similar formulas can be derived for
nonuniform Cartesian grids.

\subsection{Modification of the dimension--by--dimension WENO method}
The considerations of the previous section suggest the following modification of the
dimension--by--dimension finite volume WENO method.\\

{\bf Algorithm: Modified dimension--by--dimension WENO method}
\begin{enumerate}
\item Compute averaged values of the conserved quantities at grid cell
  interfaces using one--dimensional WENO reconstruction, i.e.\ compute
$$
Q_{i-\frac{1}{2},j}^\pm(t), \ Q_{i,j-\frac{1}{2}}^\pm(t) 
$$
at all grid cell interfaces.
\item Compute point values of the conserved quantities at the
  midpoints of grid cell
  interfaces, i.e.\ compute
$$
q_{i-\frac{1}{2},j}^\pm(t), \ q_{i,j-\frac{1}{2}}^\pm(t)
$$
using the transformation (\ref{eqn:transform}).
\item Compute fluxes at midpoints of the grid cell interfaces using a consistent
  numerical flux function, i.e.\
$$
f_{i-\frac{1}{2},j}(t) = {\cal F}(q_{i-\frac{1}{2},j}^-,
q_{i-\frac{1}{2},j}^+), \quad
g_{i,j-\frac{1}{2}}(t) = {\cal F}(q_{i,j-\frac{1}{2}}^-, q_{i,j-\frac{1}{2}}^+)
$$
\item Compute  averaged values of the flux, denoted by
  $F_{i-\frac{1}{2},j}(t)$  and $G_{i,j-\frac{1}{2}}(t)$, 
at grid cell interfaces using the transformation (\ref{eqn:transform}).  
\item Approximate the semi--discrete system (\ref{eqn:semi-discrete-fv}), using a
  high--order accurate Runge--Kutta method.
\end{enumerate}

Different versions of this method can now be considered depending 1.) on
the number of terms on the right  hand side of (\ref{eqn:transform}),
which are used to transform between cell average values and point
values and 2.) on the choice of the formula used to discretize the
derivatives in (\ref{eqn:transform}).   

We consider the following variations of the method. In the different
methods, the cell
averaged values of the conserved quantities $Q_{i-\frac{1}{2},j}^\pm$
and $Q_{i,j-\frac{1}{2}}^\pm$ are computed using either fifth or
seventh order WENO reconstruction.  
\begin{itemize}
\item {\bf method 1:} The standard dimension--by--dimension approach.
\item {\bf method 2:} Point values of the conserved quantities at midpoints of grid cell
  interfaces are computed using
\begin{equation}\label{eqn:method2-q}
\begin{split}
q_{i-\frac{1}{2},j}^\pm & = Q_{i-\frac{1}{2},j}^\pm - \frac{1}{24}
\left( Q_{i-\frac{1}{2},j-1}^\pm - 2 Q_{i-\frac{1}{2},j}^\pm +
  Q_{i-\frac{1}{2},j+1}^\pm \right)\\
q_{i,j-\frac{1}{2}}^\pm & = Q_{i,j-\frac{1}{2}}^\pm - \frac{1}{24}
\left(Q_{i-1,j-\frac{1}{2}}^\pm - 2 Q_{i,j-\frac{1}{2}}^\pm +
  Q_{i+1,j-\frac{1}{2}}^\pm \right)
\end{split}
\end{equation}
The numerical fluxes used in (\ref{eqn:semi-discrete-fv}) are computed
from point values of the flux using the relations
\begin{equation}\label{eqn:transform-flux2}
\begin{split}
F_{i-\frac{1}{2},j} & = f_{i-\frac{1}{2},j} + \frac{1}{24} \left(
f_{i-\frac{1}{2},j-1} - 2 f_{i-\frac{1}{2},j} + f_{i-\frac{1}{2},j+1}
\right)\\
G_{i,j-\frac{1}{2}} & = g_{i,j-\frac{1}{2}} + \frac{1}{24} \left(
g_{i-1,j-\frac{1}{2}} - 2 g_{i,j-\frac{1}{2}} + g_{i+1,j-\frac{1}{2}} \right)
\end{split}
\end{equation}
\item {\bf method 3:} Point values of the conserved quantities at grid
  cell interfaces are computed using 
\begin{equation}
\begin{split}
& q_{i-\frac{1}{2},j}^\pm  = Q_{i-\frac{1}{2},j}^\pm\\
 & - \frac{1}{24}
\left( - \frac{1}{8} Q_{i-\frac{1}{2},j-2}^\pm + \frac{3}{2}
  Q_{i-\frac{1}{2},j-1}^\pm - \frac{11}{4} Q_{i-\frac{1}{2},j}^\pm +
  \frac{3}{2} Q_{i-\frac{1}{2},j+1}^\pm - \frac{1}{8}
  Q_{i-\frac{1}{2},j+2}^\pm \right)\\
& - \frac{1}{1920} \left( Q_{i-\frac{1}{2},j-2}^\pm -
  4Q_{i-\frac{1}{2},j-1}^\pm + 6Q_{i-\frac{1}{2},j}^\pm - 4
  Q_{i-\frac{1}{2},j+1}^\pm + Q_{i-\frac{1}{2},j+2}^\pm \right) \\
& = Q_{i-\frac{1}{2},j}^\pm  
 - \left( -\frac{3}{640}
  Q_{i-\frac{1}{2},j-2} + \frac{29}{480} Q_{i-\frac{1}{2},j-1}^\pm -
  \frac{107}{960} Q_{i-\frac{1}{2},j}^\pm \right. \\
& \hspace*{5cm} \left. + \frac{29}{480}
  Q_{i-\frac{1}{2},j+1}^\pm- \frac{3}{640} Q_{i-\frac{1}{2},j+2}^\pm \right)
\end{split}
\end{equation}
and by an analogous formula for $q_{i,j-\frac{1}{2}}^\pm$.

The interface fluxes are computed from point values of the fluxes
using
\begin{equation}\label{eqn:transform-flux3}
\begin{split}
& F_{i-\frac{1}{2},j}  = f_{i-\frac{1}{2},j} \\
&+ \frac{1}{24} \left(
  -\frac{1}{12}f_{i-\frac{1}{2},j-2} + \frac{4}{3}
  f_{i-\frac{1}{2},j-1} - \frac{5}{2} f_{i-\frac{1}{2},j} +
  \frac{4}{3} f_{i-\frac{1}{2},j+1} -\frac{1}{12}
  f_{i-\frac{1}{2},j+2} \right)\\
& + \frac{1}{1920} \left( f_{i-\frac{1}{2},j-2} -
  4f_{i-\frac{1}{2},j-1} + 6 f_{i-\frac{1}{2},j} - 4
  f_{i-\frac{1}{2},j+1} + f_{i-\frac{1}{2},j+2} \right)
\end{split}
\end{equation}
and an analogous formula for $G_{i,j-\frac{1}{2}}$.
\end{itemize}

\begin{figure}[htb!]
 \centering 
  \begin{minipage}{0.490\linewidth}
  \centering 
    \begin{tikzpicture}[scale=0.75]
      \fill[color=white] (-3.5,-1.9) rectangle (3.5,2.9);
      \usetikzlibrary{patterns}
      \foreach \x in {-3,...,2}
      \foreach \y in {0}
      {
        \shade [shading=radial,inner color=gray!50]
        (\x,\y) rectangle (\x+1,\y+1);
      }
      \draw[very thin] (-3.5,-0.5) grid (3.5,1.5);
      \draw[dashed]  (.1,0.1) -- (0.1,0.9);
      \draw[dashed]  (-.1,0.1) -- (-0.1,0.9);
    \end{tikzpicture}
 %    \caption{Stencil to get $Q^\pm_{i+\frac{1}{2},j}$.}
  \end{minipage}
  \begin{minipage}[r]{0.490\linewidth}
  \centering     
   \begin{tikzpicture}[scale=0.75]
      \fill[color=white] (-3.5,-1.9) rectangle (3.5,2.9);
      \usetikzlibrary{patterns}
      \foreach \x in {-3,...,2}
      \foreach \y in {-1,...,1}
      {
        \shade [shading=radial,inner color=gray!50] 
        (\x,\y) rectangle (\x+1,\y+1);
      }
      \draw[very thin] (-3.5,-1.5) grid (3.5,2.5);
      \foreach \y in {-1,...,1}
      {
        \draw[dashed]  ( .1,0.1+\y)   -- ( .1,0.9+\y);
        \draw[dashed]  (-.1,0.1+\y)  -- (-.1,0.9+\y);
      }
      \fill ( .1,0.5) circle (2pt);
      \fill (-.1,0.5) circle (2pt);
        \draw[thin](0,0.5) ellipse (6pt and 4pt);
    \end{tikzpicture}
%     \caption{Stencil to get $f_{i+\frac{1}{2},j}$.}
  \end{minipage}
\caption{\label{fig:stencil1}The left plot shows the stencil for the computation of
  $Q^\pm_{i-\frac{1}{2},j}$. These averaged interface values are
  indicated by the two dashed lines. The right plot shows the stencil
  for the computation of the point values $q_{i-\frac{1}{2},j}^\pm$,
  indicated as black dots in the figure. Those point values are used
  to compute point values of the flux, denoted by
  $f_{i-\frac{1}{2},j}$. The point value of the flux is marked as
  an open ellipse.}
\end{figure}
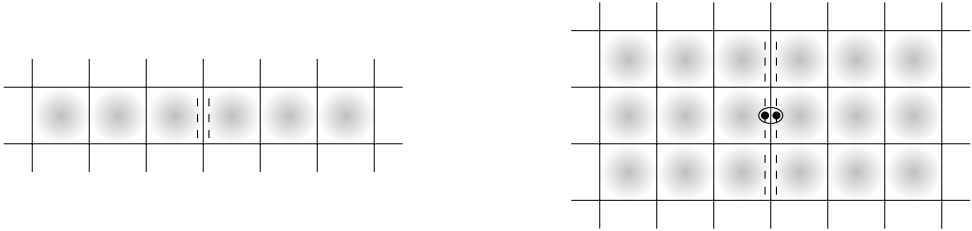

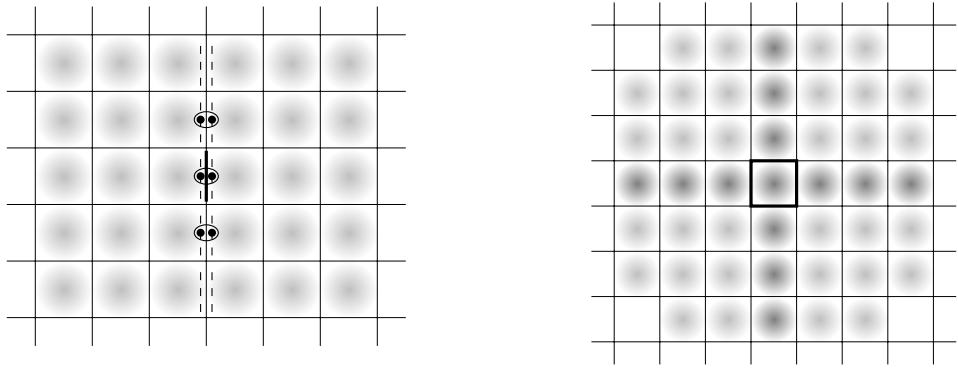
\begin{figure}[htb!]
  \begin{minipage}[r]{0.490\linewidth}
  \centering     
    \begin{tikzpicture}[scale=0.75]
      \fill[color=white] (-3.5,-2.6) rectangle (3.5,3.5);
      \usetikzlibrary{patterns}
      \foreach \x in {-3,...,2}
      \foreach \y in {-2,...,2}
      {
        \shade [shading=radial,inner color=gray!50] 
        (\x,\y) rectangle (\x+1,\y+1);
      }
      \draw[very thin] (-3.5,-2.5) grid (3.5,3.5);
      \foreach \y in {-2,...,2}
      {
        \draw[dashed]  ( .1,0.1+\y)   -- ( .1,0.9+\y);
        \draw[dashed]  (-.1,0.1+\y)  -- (-.1,0.9+\y);
      }
      \foreach \y in {-1,...,1}
      {
        \fill ( .1,0.5+\y) circle (2pt);
        \fill (-.1,0.5+\y) circle (2pt);
        \draw[thin](0,0.5+\y) ellipse (6pt and 4pt);
      }
       \draw[very thick] ( 0,0.05)   -- ( 0,0.95);
    \end{tikzpicture}
%     \caption{Stencil to get $F_{i+\frac{1}{2},j}$.}
  \end{minipage}
  \begin{minipage}{0.490\linewidth}
  \centering 
    \begin{tikzpicture}[scale=0.6]
      \fill[color=white] (-3.5,-3.6) rectangle (3.5,4.8);
      \foreach \x in {-3,...,3}
      \foreach \y in {-2,...,2}
      {
        \shade [shading=radial,inner color=gray!50]
        (\x,\y) rectangle (\x+1,\y+1);
      }
      \foreach \x in {-2,...,2}
      \foreach \y in {-3,3}
      {
        \shade [shading=radial,inner color=gray!50]
        (\x,\y) rectangle (\x+1,\y+1);
      }
      \foreach \x in {-3,...,3}
      {
        \shade [shading=radial] (0,\x) rectangle (1,\x+1);
        \shade [shading=radial] (\x,0) rectangle (\x+1,1);
      }
      \draw[very thin] (-3.5,-3.5) grid (4.5,4.5);
      \draw[very thick] (0,0) rectangle (1,1);
    \end{tikzpicture}
  \end{minipage}
\caption{\label{fig:stencil2}The left plot shows the stencil, which is used in order to
  compute the cell averaged value of the flux at the interface, i.e.\
  $F_{i-\frac{1}{2},j}$. This flux is marked as a dark solid line. The
  right plot shows the complete stencil used to update the grid cell
  in the center. The dark shaded cells are used in a
  classical dimension-by-dimension approach, i.e.\ by method 1.}
\end{figure}

In Figures \ref{fig:stencil1} and \ref{fig:stencil2} we show the
stencil used in one time stage of our method using $5$th order WENO
reconstruction together with the modification implemented in method 2.
In the left part of Figure \ref{fig:stencil1}, we show the stencil
which is used in order to compute edge averaged values of the conserved
quantities marked as two dashed lines. In a standard dimension-by-dimension approach, those edge
averaged values are used to compute the interface flux. 
In our modified method, we compute point values of the conserved
quantities using equation (\ref{eqn:method2-q}).
For this computation we need neighboring edge
averaged values, which enlarges the stencil as indicated in the right part of Figure \ref{fig:stencil1}.  
The point values of the conserved quantity (indicated by black dots)
are used to compute point values of the flux (indicated by the open
ellipse.)
%%%%%

In the left plot of Figure \ref{fig:stencil2}, we show the stencil
needed to compute edge averaged values of the flux according to
equation (\ref{eqn:transform-flux2}). This transformation requires
neighboring point values of fluxes, which further enlarges the stencil. In the
right plot of Figure \ref{fig:stencil2} we show the full stencil
of cells that are used to update one cell. The dark shaded
grid cells are those used in the classical dimension-by-dimension WENO method. 
Note that after computing all fluxes (of the modified method) for one cell, most of the work
for the neighboring cells is already done. Therefore, the larger
stencil only leads to a relatively small increase of the computational
costs, as shown below in Table \ref{table:ex2-3}. 

\begin{remark}
We can replace the flux computation in (\ref{eqn:transform-flux2}) by a
formula of the form
\begin{equation}
F_{i-\frac{1}{2},j} = f_{i-\frac{1}{2},j} + \frac{1}{24} \left(
  \bar{f}_{i-\frac{1}{2},j-1} - 2 \bar{f}_{i-\frac{1}{2},j} +
  \bar{f}_{i-\frac{1}{2},j+1} \right),
\end{equation}
and analogously for $G_{i,j-\frac{1}{2}}$, 
where $\bar{f}_{i-\frac{1}{2},k} = {\cal F}(Q_{i-\frac{1}{2},k}^-,
Q_{i-\frac{1}{2},k}^+)$, $k=j-1, j, j+1$, is a flux computed using the
averaged values of the conserved quantities. The resulting finite volume
method has a more local stencil. (For WENO5 the 16 most outer light
shaded cells in the right plot of Figure \ref{fig:stencil2} would not be
used.) However, this approach requires the
computation of two fluxes per interface. 
We have also tested such versions of the method and obtained good results. These
computations will not be presented here. 
\end{remark}

%%%%%%%%
Finally, in Table \ref{table:rates}, 
\begin{table}[htb!]
\centerline{
\begin{tabular}{|c||c|c||c|c|}
\hline
&  \multicolumn{2}{c|}{WENO5+RK5} &
 \multicolumn{2}{c|}{WENO7+RK7} \\
method & linear & nonlinear & linear & nonlinear\\
\hline
method 1 & 5 & 2 & 7 & 2 \\
method 2 & 4 & 4 & 4 & 4 \\
method 3 & 5 & 5 & 6 & 6\\
\hline
\end{tabular}}
\caption{\label{table:rates}Predicted convergence rate for the approximation of smooth problems
  with the different numerical methods.}
\end{table}
we summarize the expected convergence rates
of the different methods for the approximation of linear and nonlinear problems.

\section{Nonlinear systems in quasilinear form}\label{section:s5}
In a recent paper, Ketcheson et al.\ \cite{article:KPL2013} used WENO
reconstruction to develop high--order wave propagation methods for
hyperbolic equations in the quasilinear form
\begin{equation}\label{eqn:2d-quasilin}
q_t + A(q) q_x + B(q) q_y = 0.
\end{equation} 
Their numerical method can be written in  the semi--discrete form
\begin{equation}
\begin{split}
\frac{d}{dt} Q_{i,j}(t) = & -\frac{1}{\Delta x} \left( {\cal A}^- \Delta
  Q_{i+\frac{1}{2},j} + {\cal A}^+ \Delta Q_{i-\frac{1}{2},j} + {\cal
    A} \Delta Q_{i,j} \right)\\
&  - \frac{1}{\Delta y} \left( {\cal B}^- \Delta
  Q_{i,j+\frac{1}{2}} + {\cal B}^+ \Delta Q_{i,j-\frac{1}{2}} + {\cal
    B} \Delta Q_{i,j} \right),
\end{split}
\end{equation}
where ${\cal A}^\pm \Delta Q$ and ${\cal B}^\pm \Delta Q$ are
fluctuations which are computed using an eigenvector decomposition
of the jump of the piecewise polynomial reconstructed quantity $q$ at
each grid cell interface as explained in \cite{article:KPL2013,book:RJL2002}.
The two--dimen\-sion\-al version of the method in \cite{article:KPL2013}
is based on a dimension--by--dimension WENO reconstruction of $q$. It
is second order accurate for smooth solutions of nonlinear hyperbolic
systems in analogy to the dimension--by--dimension approach for
nonlinear hyperbolic systems in divergence form.  

The approach suggested in this paper to obtain high order of accuracy 
can be extended to hyperbolic systems of the form
(\ref{eqn:2d-quasilin}).
The computation of the terms ${\cal A}^\pm \Delta Q$ and ${\cal B}^\pm
\Delta Q$ can be done in analogy to the flux computation, i.e.\ by
first computing high order point values of $q$ at the interfaces,
by using these values to compute high--order accurate point values of the
fluctuations, and finally by computing grid cell interface averaged
values of the fluctuations from the point values of the fluctuations.  

The discretization of ${\cal A}\Delta
Q_{i,j}$ and ${\cal B}\Delta Q_{i,j}$ requires  some additional transformations.
Consider the discretization of
\begin{equation} \label{eqn:2dint}
{\cal A}\Delta Q_{i,j}  \approx \frac{1}{\Delta y}
\int_{y_{j-\frac{1}{2}}}^{y_{j+\frac{1}{2}}}
\int_{x_{i-\frac{1}{2}}}^{x_{i+\frac{1}{2}}} A(q(x,y,t)) q_x(x,y,t) dx
dy.
\end{equation} 
Let $x_{i-\frac{1}{2}} \le x_1 < \ldots < x_\ell \le
x_{i+\frac{1}{2}}$ and $c_1,
\ldots, c_\ell$ denote the nodes and weights of a quadrature formula,
which can be written in the general form
\begin{equation}\label{eqn:qf}
Q[f] = \sum_{k=1}^\ell c_k f(x_k) \approx
\int_{x_{i-\frac{1}{2}}}^{x_{i+\frac{1}{2}}} f(x) dx.
\end{equation}
Furthermore,
let $q_{i,j}^1(x)$, $x_{i-\frac{1}{2}}< x < x_{i+\frac{1}{2}}$ denote
the $p$-th oder accurate WENO
reconstruction in the $x$--direction of the quantity $q$ in cell
$(i,j)$, compare with Section \ref{section:s2}. The evaluation of
$q_{i,j}^1$ at a quadrature node provides us with a point value of $q$ in the
$x$-direction and an averaged value of $q$
in the $y$--direction, i.e.\
\begin{equation}
q_{i,j}^1(x_k) = \frac{1}{\Delta y}
\int_{y_{j-\frac{1}{2}}}^{y_{j+\frac{1}{2}}} q(x_k,y) dy + {\cal
  O}(\Delta x^p + \Delta y^p) \quad k=1,\ldots,\ell.
\end{equation}
Furthermore, differentiating the polynomial $q^1_{i,j}$ provides us with an
approximation of the averaged derivative, i.e.\
\begin{equation}
\left( q^1_{i,j} \right)_x (x_k) = \frac{1}{\Delta y}
\int_{y_{j-\frac{1}{2}}}^{y_{j+\frac{1}{2}}} q_x(x_k, y) dy + {\cal
  O}(\Delta x^{p-1} + \Delta y^p), \quad k=1,\ldots,\ell.
\end{equation}

Using the transformation from averaged values to point values, i.e.\
(\ref{eqn:transform}), we compute approximations of the point values
$q(x_k,y_j)$ and $q_x(x_k,y_j)$ for $k=1, \ldots, \ell$ and $y_j = (y_{j-\frac{1}{2}}+y_{j+\frac{1}{2}})/2$. 
These point values are computed using neighboring
averaged values, i.e.\ $q_{i,j-2}^1(x_k)$, $q_{i,j-1}^1(x_k)$, $q_{i,j+1}^1(x_k)$
and $q_{i,j+2}^1(x_k)$ for $k=1, \ldots, \ell$ and analogously for the
derivatives.
Now we can evaluate the point values
$A(q_{i,j}^1(x_k,y_j,t))q^1_x(x_k,y_j,t)$ for $k=1, \ldots, \ell$. 
Using again the transformation
(\ref{eqn:transform}), we compute averaged values of these quantities in
the $y$--direction and denote them by
\begin{equation}
\overline{A(q_{i,j}^1(x_k,t)) q^1_x(x_k)} \approx \frac{1}{\Delta y}
\int_{y_{j-\frac{1}{2}}}^{y_{j+\frac{1}{2}}} A(q(x_k,y)) q_x(x_k,y) dy
\quad k=1, \ldots, \ell.
\end{equation}  
These values are finally used in a one--dimensional quadrature formula of the form
(\ref{eqn:qf}), giving an approximation of the two--dimensional 
integral in (\ref{eqn:2dint}):
\begin{equation}
\sum_{k=1}^\ell c_k \overline{A(q_{i,j}^1(x_k,t)) q^1_x(x_k)} =
{\cal A} \Delta Q_{i,j}.
\end{equation}
Analogously we can compute the term 
\begin{equation}
{\cal B} \Delta Q_{i,j} \approx \frac{1}{\Delta x}
\int_{y_{j-\frac{1}{2}}}^{y_{j+\frac{1}{2}}}
\int_{x_{i-\frac{1}{2}}}^{x_{i+\frac{1}{2}}} B(q(x,y,t)) q_y(x,y,t) dx dy.
\end{equation} 

\section{Extension to higher dimensions}
\label{sec:extens-high-dimens}
We now present the basic formula which is needed 
 for an extension of the modified WENO method to
higher dimensions. Let us first introduce some additional notation.
For $\boldsymbol{x} \in \mathbb{R}^d,\boldsymbol{n} \in \mathbb{N}^d,$
let
$\boldsymbol{x}^{\boldsymbol{n}}= x_1^{n_1}\cdots x_d^{n_d},
|\boldsymbol{n}|=n_1+\cdots +n_d$,
$\boldsymbol{n}!=n_1!\cdots n_d!$ and
$\boldsymbol{n+1}! = (n_1+1)! \cdots (n_d+1)!$.
Furthermore we restrict our
considerations to equidistant Cartesian grids and set
$h=\Delta x_1=...=\Delta x_d$.
For
$\boldsymbol{x}_{\boldsymbol{i}},\boldsymbol{\xi} \in \mathbb{R}^d$
with $|\xi_j|\leq h,$ $ j=1,...,d$,
the multidimensional Taylor series expansion yields
\begin{equation}
 \label{eq:multi-taylor}
 \begin{aligned}
 q(\boldsymbol{x}_{\boldsymbol{i}}+\boldsymbol{\xi})
 &= \sum_{\substack{0\leq n_1,...,n_d, \\n_1+\cdots +n_d\leq p}}
 \frac{\xi_1^{n_1}\cdots\xi_d^{n_d}}{n_1!\cdots n_d!}
 \frac{\partial^{n_1+\cdots +n_d}}
       {\partial x^{n_1}_{1}\cdots\partial x^{n_d}_{d}}
 q(\boldsymbol{x}_{\boldsymbol{i}}) + \mathcal{O}(h^{p+1})\\
 &= \sum_{\substack{0\leq  |\boldsymbol{n}| \leq p}}
 \frac{\boldsymbol{\xi}^{\boldsymbol{n}}}{\boldsymbol{n}!}
 \frac{\partial^{|\boldsymbol{n}|}}
       {\partial\boldsymbol{x}^{\boldsymbol{n}}}
 q(\boldsymbol{x}_{\boldsymbol{i}}) + \mathcal{O}(h^{p+1}).
 \end{aligned}
\end{equation}
For $\xi \in \mathbb{R}, n \in \mathbb{N}$ we get
\begin{equation}
 \label{eq:1}
 \int_{ -\frac{h}{2}}^{ \frac{h}{2}} \xi^{n}d\xi
 =\left.\frac{1}{n+1}\xi^{n+1}\right|_{ -\frac{h}{2}}^{ \frac{h}{2}}
 =
 \begin{cases}
   0 &\text{ if } n \text{ is odd},\\
   \frac{h^{n+1}}{2^n (n+1)} &\text{ if } n \text{ is even}.
 \end{cases}
\end{equation}
For an even number $p$ this leads to
\begin{equation}
 \label{eq:multi-d}
 \begin{aligned}
   Q_{\boldsymbol{i}}
   &=\frac{1}{h^d}\int_{x_{i_1} -\frac{h}{2}}^{x_{i_1} +\frac{h}{2}}
   \cdots\int_{x_{i_d} -\frac{h}{2}}^{x_{i_d} +\frac{h}{2}}
   q(x_1,...,x_d) dx_d...dx_1\\
   &=\frac{1}{h^d}\int_{-\frac{h}{2}}^{ \frac{h}{2}}
   \cdots\int_{ -\frac{h}{2}}^{ \frac{h}{2}}
   q(x_{i_1}+\xi_1,...,x_{i_d}+\xi_d) d\xi_d...d\xi_1\\
   &=\frac{1}{h^d}
   \sum_{\substack{0\leq|\boldsymbol{n}| \leq p}}
   \frac{\partial^{|\boldsymbol{n}|}}
       {\partial \boldsymbol{x}^{\boldsymbol{n}}}
   q(\boldsymbol{x}_{\boldsymbol{i}})
   \int_{-\frac{h}{2}}^{ \frac{h}{2}}
   \cdots\int_{ -\frac{h}{2}}^{ \frac{h}{2}}
\frac{\xi_1^{n_1}\cdots\xi_d^{n_d}}{\boldsymbol{n}!}d\xi_d...d\xi_1
   + \mathcal{O}(h^{p+1})    \\
   &=\frac{1}{h^d}\sum_{\substack{0\leq|\boldsymbol{n}| \leq p,\\ n_j=2k_j}}
   \frac{h^{|\boldsymbol{n}|+d}}{2^{|\boldsymbol{n}|} {(\boldsymbol{n+1})!}}
   \frac{\partial^{|\boldsymbol{n}|}}
       {\partial \boldsymbol{x}^{\boldsymbol{n}}}
   q(\boldsymbol{x}_{\boldsymbol{i}})
   + \mathcal{O}(h^{p+2})    \\
   &=q(\boldsymbol{x}_{\boldsymbol{i}})
   +\sum_{\substack{0<|\boldsymbol{n}| \leq p,\\ n_j=2k_j}}
   \frac{h^{|\boldsymbol{n}|}}{2^{|\boldsymbol{n}|} {(\boldsymbol{n+1})!}}
   \frac{\partial^{|\boldsymbol{n}|}}
       {\partial \boldsymbol{x}^{\boldsymbol{n}}}
   q(\boldsymbol{x}_{\boldsymbol{i}})
   + \mathcal{O}(h^{p+2}).
 \end{aligned}
\end{equation}
By approximating the corresponding derivatives, we can now obtain
transformation formulas for any dimension and order.
For $p=2$ we retain method 2 and for $p=4$ we retain method 3.
But notice, while method 2 is a simple sum of second derivatives,
higher order transformations contain also cross terms when applied in
more then one dimensions. Note that at grid cell interfaces of a two-dimensional
Cartesian mesh, we need to apply the one-dimensional transformation formulas
($d=1$)
and at grid cell interfaces of a three-dimensional cell we apply the
two-dimensional formulas ($d=2$).

\section{Numerical results for the Euler equations of gas dynamics}
\label{section:results}
In this section we present different simulations and convergence
studies. We use the two--dimensional Euler equations of gas dynamics
as our model problem, i.e.\ we consider approximations of
\begin{equation}\label{euler2d}
\partial_t \left( \begin{array}{c}
\rho\\ \rho u\\ \rho v \\ E\end{array}\right) + \partial_x \left(
\begin{array}{c}
\rho u \\ \rho u^2 + p\\ \rho u v \\ u(E+p)\end{array}\right)
+ \partial_y \left( \begin{array}{c}
 \rho v \\ \rho u v \\ \rho v^2 + p\\ v(E+p)\end{array}\right)
= 0,
\end{equation}
with the ideal gas equation of state
\begin{equation}
E = \frac{p}{\gamma-1} + \frac{1}{2} \rho (u^2 + v^2).
\end{equation}
The initial values will be specified below for each test problem. We
always set $\gamma = 1.4$.

In the following subsections, we will present different tabels with convergence
studies. There, the $\| \cdot \|_1$-norm  of the error in density is
shown for different grids. If an exact solution is available, we use
those as reference solution. Otherwise, the reference solution is a
numerical solution computed on a highly refined mesh. We compute the
experimental order of convergence using the formula
\begin{equation}
EOC = \frac{\log (\|\rho_m - \rho_{ref} \| / \| \rho_{2m} -
  \rho_{ref}\|)}{\log 2},
\end{equation}
where the index $m$ indicates the number of grid points in the $x$ and
the $y$ direction. 
Note that for the other conserved variables we alway obtained comparable results which
are not shown here. In all computations we used time steps corresponding to 
$CFL \approx 0.9$.

\subsection{Smooth test problems and convergence studies}

\subsubsection{Linear problem}
\begin{example}\label{example:1} We consider periodic solutions of
  (\ref{euler2d}) on the domain  $[0,1]\times[0,1]$. 
  The initial values are given by 
  \begin{equation}
    \begin{split}
      \rho(x,y,0) & = 1 + 0.5 \sin(2 \pi x) \cos(2 \pi y) \\
      p(x,y,0) & = 1 \\
      u(x,y,0) & = v(x,y,0) = 1.
    \end{split}
  \end{equation}
  In this case, velocity and pressure remain constant for all times
  and the initial density profile is advected by the velocity field.  
  Thus we are approximating a problem in the linear regime. 
\end{example}
In Table \ref{table:ex1} we show results of a numerical convergence
study. Here we compute the $\| \cdot \|_1$--norm of the error in
density by comparing the solution obtained on different grids with the 
exact solution. 
In Table \ref{table:ex1} we show results with fifth order WENO-Z
reconstruction, using $\epsilon = \Delta x^4$ and $p=2$, compare with
Appendix \ref{appendix:s2}. 
In all of these computations, we used RK5 as time stepping scheme, see
Appendix \ref{app:RK4-6}.

As expected in the linear case, the simple dimension-by-dimension
approach, implemented in method 1,  converges with fifth order. 
The full order of convergence of the WENO-Z reconstruction is also
retained by method 3. 
By using method 2, we observe a loss of accuracy and, as expected, a convergence
rate of four.
\begin{table}[htb]
  \centerline{
    \begin{tabular}{|c||c|c|c|c|c|c|}
      \hline
      & \multicolumn{2}{c|}{method 1} 
      & \multicolumn{2}{c|}{method 2}
      & \multicolumn{2}{c|}{method 3}\\
      grid & $\|\rho-\rho_{exact}\|_1$ & EOC & $\|\rho-\rho_{exact}\|_1$ &
      EOC & $\|\rho-\rho_{exact}\|_1$ & EOC\\ \hline
      $64^2$  & 4.61953d-7 & &  5.65033d-7 & & 4.61965d-7 &\\
      $128^2$ & 1.44675d-8 & {\bf 5.00 }  & 2.49135d-8 & {\bf 4.50  } &
      1.44677d-8 & {\bf 5.00 } \\
      $256^2$ & 4.52368d-10 & {\bf 5.00 } & 1.34356d-9  & {\bf 4.21  }&
      4.52369d-10& {\bf 5.00 } \\
      $512^2$ & 1.41384d-11&  {\bf 5.00 } & 8.01169d-11 & {\bf  4.07 }&
      1.41384d-11 & {\bf 5.00 } \\
      \hline
    \end{tabular}}
  \caption{\label{table:ex1}
    Convergence study for Example \ref{example:1} with $5$th order
    WENO-Z reconstruction and RK5.
    For these computations we used time steps with $CFL\approx0.9$ and
    the Lax-Friedrichs flux.}
\end{table} 
In Table \ref{table:ex1-2} we show the results where the seventh order
WENO-Z was combined with RK7, see Appendix
\ref{appendix:s2} and \ref{app:RK4-6}.  
\begin{table}[htb]
  \centerline{
    \begin{tabular}{|c||c|c|c|c|c|c|}
      \hline
      & \multicolumn{2}{c|}{method 1} 
      & \multicolumn{2}{c|}{method 2}
      & \multicolumn{2}{c|}{method 3}\\
      grid & $\|\rho-\rho_{exact}\|_1$ & EOC 
      & $\|\rho-\rho_{exact}\|_1$ &
      EOC & $\|\rho-\rho_{exact}\|_1$ & EOC\\ \hline
      $64^2$  & 9.52785d-10 & &  3.21588d-7 & & 1.11271d-9 &\\
      $128^2$ & 7.46432d-12 &{\bf 7.00}& 2.01439d-08 &{\bf 4.00}& 
      1.16519d-11 & {\bf 6.58 } \\
      $256^2$ & 5.84576d-14 &{\bf 7.00}& 1.25974d-09 &{\bf 4.00}&
      1.51789d-13 & {\bf 6.26 } \\
      $512^2$ & 4.56021d-16 &{\bf 7.00}& 7.87454d-11 &{\bf 4.00}&
      2.22517d-15 & {\bf 6.09 } \\
      \hline
    \end{tabular}}
  \caption{\label{table:ex1-2}
    Convergence study for Example \ref{example:1} with $7$th order 
    WENO-Z reconstruction and RK7.
    For these computations we used time steps with $CFL\approx0.9$ and
    the Lax-Friedrichs flux.}
\end{table}
Again we observe that the order of convergence is as expected for each
of the methods. 
Note that for the combinations {\em WENO-Z5+RK5+method 2} and
{\em WENO-Z7+RK7+method 3} we observe, at least on coarser grids, an experimental order of
convergence which is above the theoretically expected order of convergence.

\subsubsection{Nonlinear problems}
\begin{example}\label{example:2}
Now we consider the two-dimensional vortex evolution problem (see e.g.\
\cite{article:HS1999}) on the periodic domain $[-7,7]\times[-7,7].$
The initial data consist of a mean flow $\rho=u=v=p=1$, which is perturbed
by adding
\begin{equation}
  \begin{pmatrix}
    \delta \rho \\\delta u \\\delta v \\\delta p 
  \end{pmatrix}
  =
  \begin{pmatrix}
    (1+\delta T)^{1/(\gamma-1)}-1 \\
    -y\frac{\sigma}{2\pi}e^{0.5(1-r^{²})}\\
    x\frac{\sigma}{2\pi}e^{0.5(1-r^{²})}\\
    (1+\delta T)^{\gamma/(\gamma -1)}-1
  \end{pmatrix}.
\end{equation}
Here $\delta T$, the perturbation in the temperature, is given by
\begin{equation}
  \delta T =-\frac{(\gamma-1)\sigma^2}{8\gamma\pi^2}e^{1-r^2},  
\end{equation}
with $r^2 = x^2 + y^2$ and the vortex strength $\sigma = 5$.
\end{example}
In the Tables \ref{table:ex2-1} and \ref{table:ex2-2}, we show the error
and the experimental convergence rates for the approximation of
density using the three different methods.
Here the exact solution, which at time $t=14$ agrees with the initial
values, was used as reference solution.

In Table \ref{table:ex2-1}, we use the fifth
order accurate WENO-Z reconstruction together with RK5, and
in Table \ref{table:ex2-2} we use the
seventh order accurate WENO-Z reconstruction with RK7.    
With the dimension--by--dimension approach implemented in method 1, 
we obtain similar results for both reconstructions. 
By refining the grid, the experimental order of convergence (EOC)
drops to second order, as expected for nonlinear problems.
Method 2 and method 3 give almost identical results on coarser
grids for the fifth order WENO-Z reconstruction.
In both cases we obtain a convergence rate of about five. 
Only on very fine grids the lower order of method 2 became visible. 
By using seventh order WENO-Z reconstruction, the drop in the
convergence rate of method 2 is much more obvious and can already be
seen on coarser grids.  
As already seen in Example \ref{example:1}, the combinations
{\em WENO-Z5+RK5+method 2} and {\em WENO-Z7+RK7+method 3} provide even
better results than expected.

\begin{table}[htb]
  \centerline{
    \begin{tabular}{|c||c|c|c|c|c|c|}
      \hline
      & \multicolumn{2}{c|}{method 1} 
      & \multicolumn{2}{c|}{method 2} 
      & \multicolumn{2}{c|}{method 3}\\
      grid & $\|\rho-\rho_{exact}\|_1$ & EOC 
           & $\|\rho-\rho_{exact}\|_1$ & EOC 
           & $\|\rho-\rho_{exact}\|_1$ & EOC \\ \hline
      $64^2$  & 2.07461d-4 &          & 1.55933d-4 &    &
      1.55288d-4 & \\
      $128^2$ & 2.95314d-5 &{\bf 2.81}& 8.17206d-6 &{\bf 4.25  } & 
      8.15400d-6 & {\bf 4.25 } \\
      $256^2$ & 7.03771d-6 &{\bf 2.07}& 2.40376d-7 &{\bf 5.09  } & 
      2.36830d-7 & {\bf 5.11 } \\
      $512^2$ & 1.75592d-6 &{\bf 2.00}& 7.71743d-9 &{\bf 4.96 }  &
      7.40743d-9 & {\bf 5.00  } \\
      $1024^2$& 4.39556d-7 &{\bf 2.00}& 2.57207d-10&{\bf 4.91  } &
      2.30667d-10& {\bf 5.01  } \\
      $2048^2$& 1.09902d-7 &{\bf 2.00}& 9.41754d-12&{\bf 4.77  } &
      7.20526d-12& {\bf 5.00  } \\
      $4096^2$& 2.74756d-8 &{\bf 2.00}& 4.00267d-13&{\bf 4.56  } &
      2.27363d-13& {\bf 4.99  } \\
      \hline
    \end{tabular}}
  \caption{\label{table:ex2-1}
  Convergence study for Example \ref{example:2} with fifth order
  accurate WENO-Z reconstruction and fifth order accurate Runge-Kutta 
  method.}
\end{table} 
\begin{table}[htb]
  \centerline{
    \begin{tabular}{|c||c|c|c|c|c|c|}
      \hline
      & \multicolumn{2}{c|}{method 1} 
      & \multicolumn{2}{c|}{method 2} 
      & \multicolumn{2}{c|}{method 3}\\
      grid & $\|\rho-\rho_{exact}\|_1$ & EOC 
           & $\|\rho-\rho_{exact}\|_1$ & EOC 
           & $\|\rho-\rho_{exact}\|_1$ & EOC \\ \hline
      $64^2$  & 1.18318d-4 &          & 3.18576d-5 &          & 
      3.09485d-5 & \\
      $128^2$ & 2.80558d-5 &{\bf 2.08}& 6.08063d-7 &{\bf 5.71}&
      4.85418d-7 & {\bf  5.99 } \\
      $256^2$ & 7.01366d-6 &{\bf 2.00}& 1.96245d-8 &{\bf 4.95}&
      4.34579d-9 & {\bf  6.80 } \\
      $512^2$ & 1.75733d-6 &{\bf 2.00}& 1.20938d-9 &{\bf 4.02}&
      3.46585d-11 & {\bf 6.97  } \\
      $1024^2$& 4.39593d-8 &{\bf 2.00}& 7.57974d-11&{\bf 4.00}&
      2.77524d-13 & {\bf 6.96  } \\
      $2048^2$& 1.10050d-8 &{\bf 2.00}& 4.73921d-12&{\bf 4.00}&
      2.40689d-15 & {\bf 6.85  } \\
      \hline
\end{tabular}}
\caption{\label{table:ex2-2}
  Convergence study for Example \ref{example:2} with seventh order
  accurate WENO-Z reconstruction and seventh order accurate
  Runge-Kutta method.}
\end{table} 
In Table \ref{table:ex2-3} we show the result of a performance-test.
In order to keep the test as fair as possible we ran all
performance-tests on the same machine and always exclusively.
The code has been parallelized using OpenMP and each test was
performed with two threads. 
We see that in our implementation the additional computations, needed
by method 2 and method 3, increase the computational costs on average
by 10\%--16\%. 
Furthermore, we observe that the relative influence is smaller for the
more expensive seventh order WENO-Z reconstruction.
This makes sense, since the
additional computations required by our modified WENO methods 
are independent of the chosen reconstruction. 
\begin{table}[htb]
  \centerline{
    \begin{tabular}{|c||c|c|c||c|c|c|}
      \hline
      & \multicolumn{3}{c||}{WENO-Z5 + RK5} 
      & \multicolumn{3}{c| }{WENO-Z7 + RK7}\\
      grid & method 1 &  method 2  & method 3 
           & method 1 &  method 2  & method 3  \\ \hline
      $64^2$  & 1.00 & 1.09 & 1.12 & 1.00 & 1.07 & 1.09 \\
      $128^2$ & 1.00 & 1.15 & 1.18 & 1.00 & 1.10 & 1.12 \\
      $256^2$ & 1.00 & 1.15 & 1.17 & 1.00 & 1.10 & 1.12 \\
      $512^2$ & 1.00 & 1.15 & 1.17 & 1.00 & 1.10 & 1.12 \\
      \hline
      $ \varnothing   $ & 1.00 & 1.14 & 1.16 & 1.00 & 1.10 & 1.11 \\
      \hline
\end{tabular}}
\caption{\label{table:ex2-3}
  Performance-test on a Intel Core 2 Duo CPU E6850 with 3.00GHz and
  2GB RAM for Example \ref{example:2}. Run time normalized with respect
  to method 1 on each grid. Compare with a performance test of the
  method studied by Zhang, Zhang, Shu in \cite{article:ZZS2011}.}
\end{table} 

\begin{figure}
\includegraphics[width=0.45\textwidth]{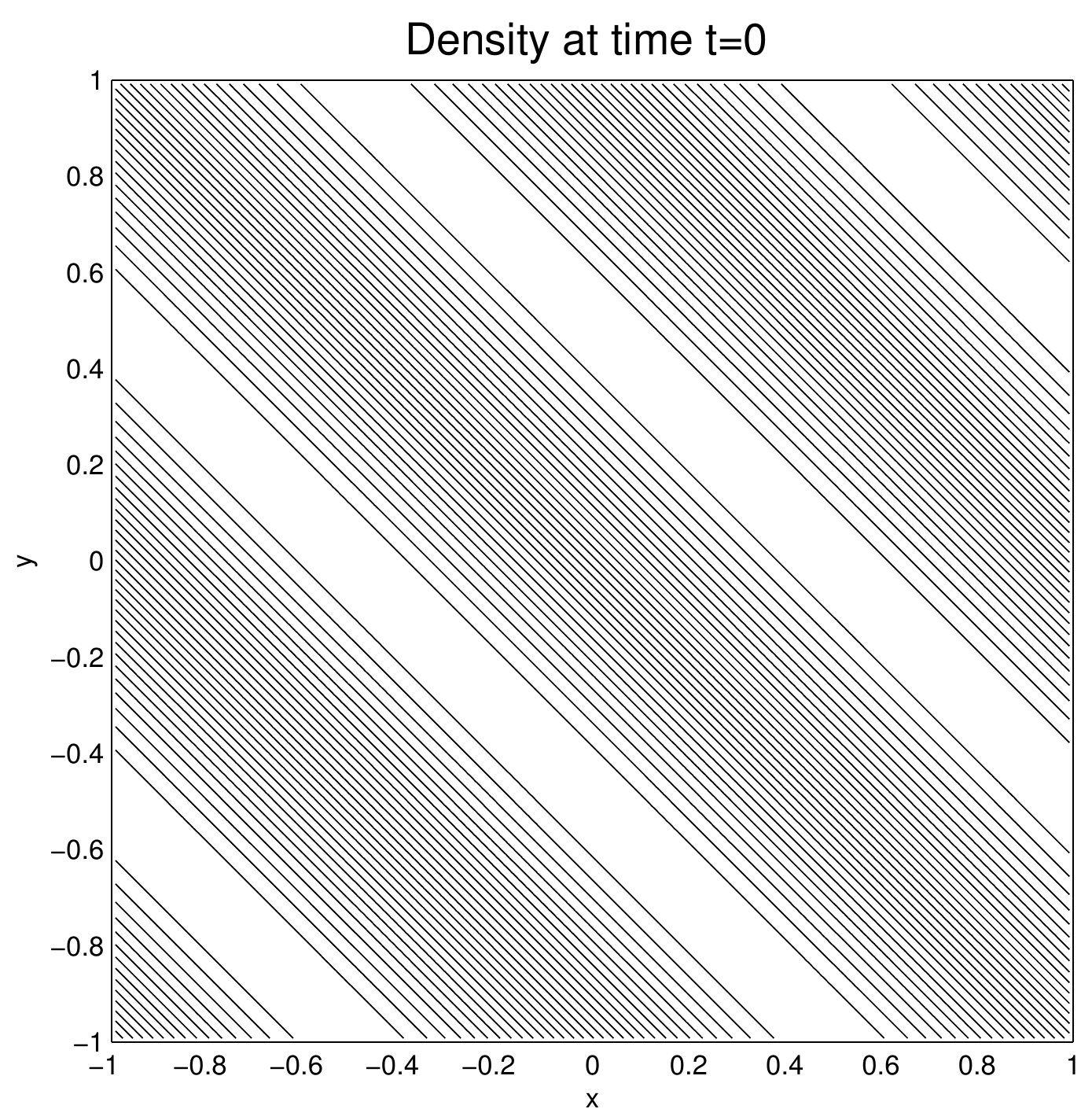}\hfill
\includegraphics[width=0.45\textwidth]{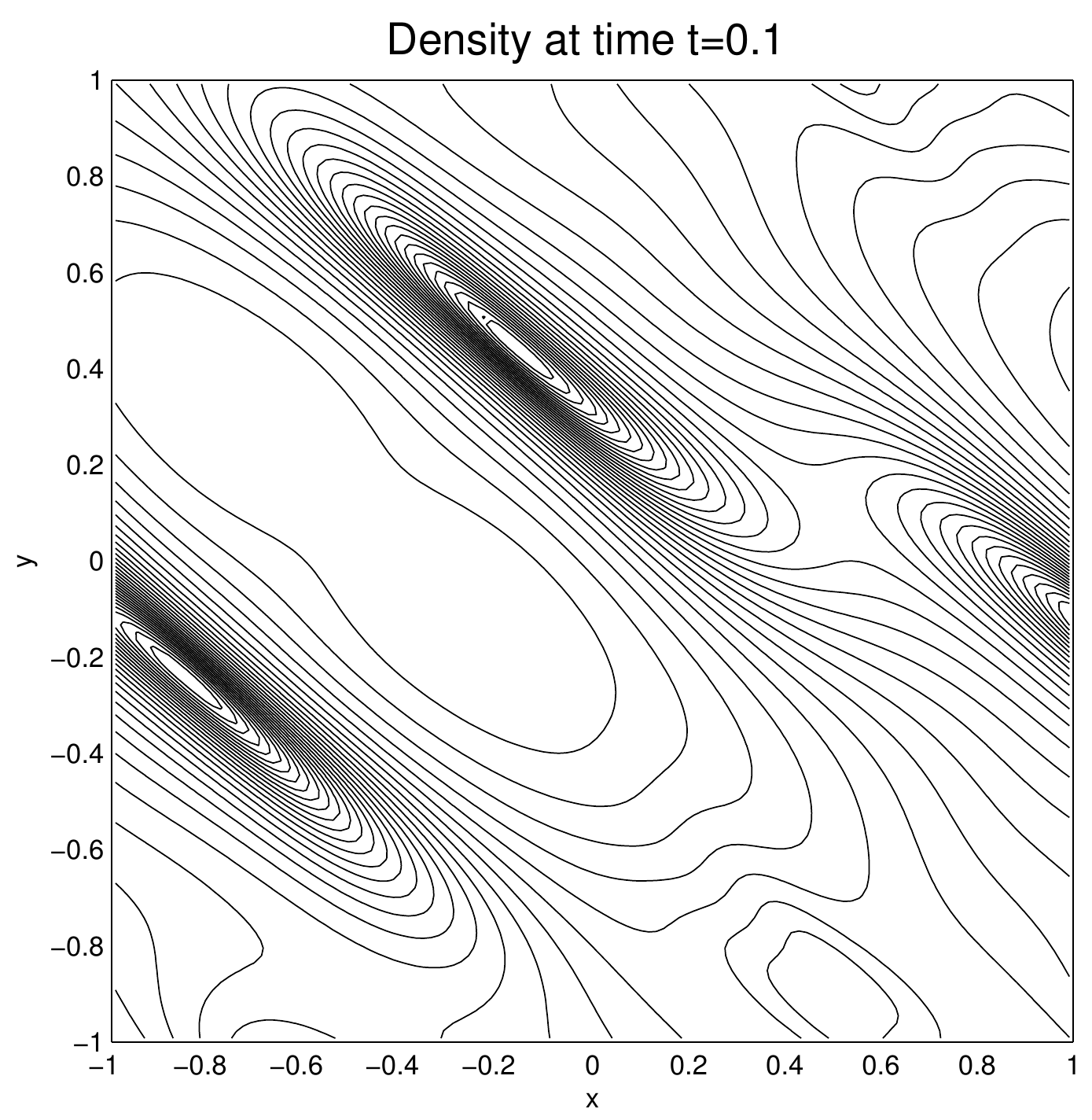}
\caption{\label{figure:ex2b}
  Contour plots for initial data and numerical
  solution at time $t=0.1$ for Example \ref{example:2b} using
  $128 \times 128$ grid cells.}   
\end{figure}

\begin{example}\label{example:2b}
Now we solve the Euler equations with initial values of the form
\begin{equation}
\begin{split}
\rho(x,y,0) & = 1 + \frac{1}{2} \sin(\pi(x+y))\\
u(x,y,0) & = \cos(\pi(x+2y))\\
v(x,y,0) & = 1-\frac{1}{2}\sin(\pi(2x+y))\\
p(x,y,0) & = 1+\frac{1}{2}\sin(\pi(x-y))
\end{split}
\end{equation}
in the domain $[-1,1]\times[-1,1]$ with periodicity condition. The
solution will be computed at time $t=0.1$ and compared to the solution at
different grids with a reference solution computed using method 3 with
seventh order reconstruction and RK7 on a grid with $4096 \times 4096$
grid cells.
\end{example}
Figure \ref{figure:ex2b} shows the solution structure for Example \ref{example:2b}.
In Table \ref{table:ex2b-1} we present a convergence study for Example
\ref{example:2b} using the three different methods with WENO-Z5 and
RK5.  
We observe the same behavior as in Example \ref{example:2}.
With the standard dimension--by--dimension approach, the method
converges with second order.
We clearly see the improved accuracy of method 2 and
method 3 compared to the standard approach used in method 1. 
But only on the very finest grid we see a considerable difference
between method 2 and method 3.
\begin{table}[htb!]
  \centerline{
    \begin{tabular}{|c||c|c|c|c|c|c|}
      \hline
      & \multicolumn{2}{c|}{method 1} 
      & \multicolumn{2}{c|}{method 2} 
      & \multicolumn{2}{c|}{method 3}\\
      grid & $\|\rho-\rho_{ref}\|_1$ & EOC 
           & $\|\rho-\rho_{ref}\|_1$ & EOC 
           & $\|\rho-\rho_{ref}\|_1$ & EOC \\ \hline
      $64^2$  & 2.00097d-3 &          & 1.32954d-3 &          &
      1.36177d-3 & \\
      $128^2$ & 4.88291d-4 &{\bf 2.03}& 7.77433d-5 &{\bf 4.10}&
      8.06215d-5 &{\bf 4.08} \\
      $256^2$ & 1.25512d-4 &{\bf 1.96}& 2.57926d-6 &{\bf 4.91}& 
      2.76102d-6 &{\bf 4.87} \\
      $512^2$ & 3.15995d-5 &{\bf 1.99}& 7.48302d-8 &{\bf 5.11}&
      8.43834d-8 &{\bf 5.03} \\
      $1024^2$& 7.90940d-6 &{\bf 2.00}& 2.25639d-9 &{\bf 5.05}&
      2.55045d-9 &{\bf 5.05} \\
      $2048^2$& 1.97774d-6 &{\bf 2.00}& 9.55098d-11&{\bf 4.56}&
      7.72001d-11&{\bf 5.05} \\
      \hline
    \end{tabular}}
  \caption{\label{table:ex2b-1}
    Convergence study for Example \ref{example:2b} with fifth order
    accurate WENO-Z reconstruction and fifth order accurate
    Runge-Kutta method.  
    The reference solution was computed using method 3 with seventh
    order WENO-Z reconstruction and RK7 on a grid with 
    $4096 \times 4096$ grid cells.}
\end{table} 

In Table \ref{table:ex2b-2} we show the results of a numerical
convergence study of Example \ref{example:2b} using the three
different methods with WENO-Z7 and RK7.
For the dimension--by--dimension approach, the higher order spatial
reconstruction did not lead to any increase in accuracy, compared to
method 1 with WENO-Z5. 
The same holds for method 2 on the finest grids.
For method 2 on coarse grids and for method 3 we do observe the gain
in accuracy due to the higher order reconstruction.
With this higher order spatial reconstruction we also observe the
expected fourth order convergence rate of method 2. 
Furthermore, we observe that the use of method 3 leads to a smaller
error if compared to method 2.   

\begin{table}[htb!]
  \centerline{
    \begin{tabular}{|c||c|c|c|c|c|c|}
      \hline
      & \multicolumn{2}{c|}{method 1} 
      & \multicolumn{2}{c|}{method 2} 
      & \multicolumn{2}{c|}{method 3}\\
      grid & $\|\rho-\rho_{ref}\|_1$ & EOC 
           & $\|\rho-\rho_{ref}\|_1$ & EOC 
           & $\|\rho-\rho_{ref}\|_1$ & EOC \\ \hline
      $64^2$  & 1.84119d-3 &          & 5.01165d-4 &          &
      5.15434d-4 & \\
      $128^2$ & 5.01723d-4 &{\bf 1.88}& 1.14852d-5 &{\bf 5.45}&
      1.01807d-5 &{\bf 5.66} \\
      $256^2$ & 1.26335d-4 &{\bf 1.99}& 4.44425d-7 &{\bf 4.69}&
      1.06754d-7 &{\bf 6.58} \\
      $512^2$ & 3.16297d-5 &{\bf 2.00}& 2.85272d-8 &{\bf 3.96}&
      8.64458d-10&{\bf 6.95} \\
      $1024^2$& 7.91038d-6 &{\bf 2.00}& 1.79754d-9 &{\bf 3.99}&
      6.60027d-12 &{\bf 7.03} \\
      $2048^2$& 1.97777d-6 &{\bf 2.00}& 1.12503d-10&{\bf 4.00}&
      5.25443d-14 &{\bf 6.97} \\
      \hline
    \end{tabular}}
  \caption{\label{table:ex2b-2}
    Convergence study for Example \ref{example:2b} with seventh order 
    accurate WENO-Z reconstruction and seventh order accurate
    Runge-Kutta method. 
    The reference solution was computed using method 3 on a grid with
    $4096 \times 4096$ grid cells.}
\end{table}

\subsection{Test problems with discontinuous solutions}
We tested our methods for the standard two-dimensional Riemann
problems proposed by Schultz-Rinne \cite{article:Schultz-Rinne1993} and
did not observe any numerical problems by using the proposed
methods for problems with discontinuities.

\begin{figure}[htb!]
\includegraphics[width=0.45\textwidth]{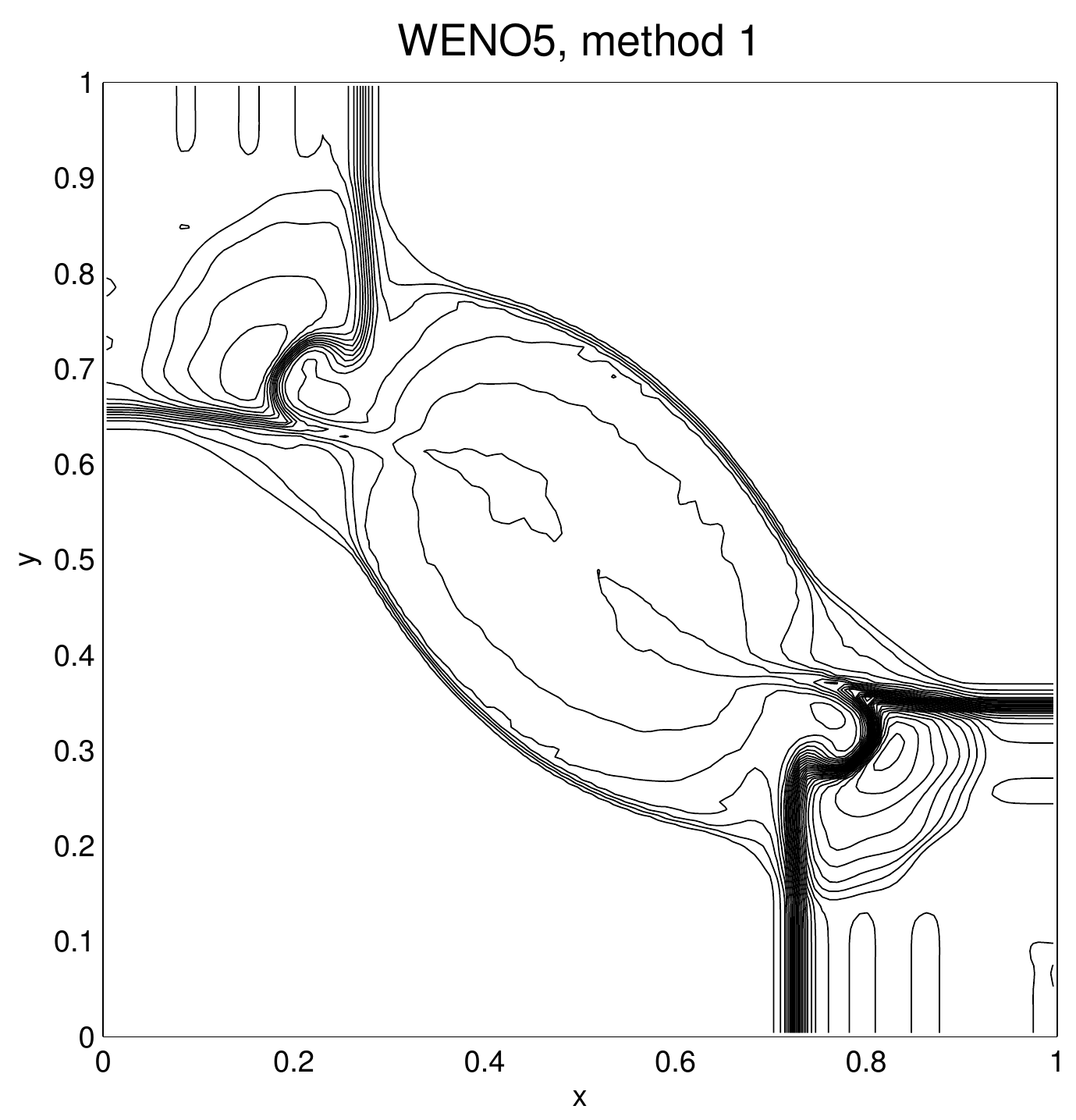}\hfil
\includegraphics[width=0.45\textwidth]{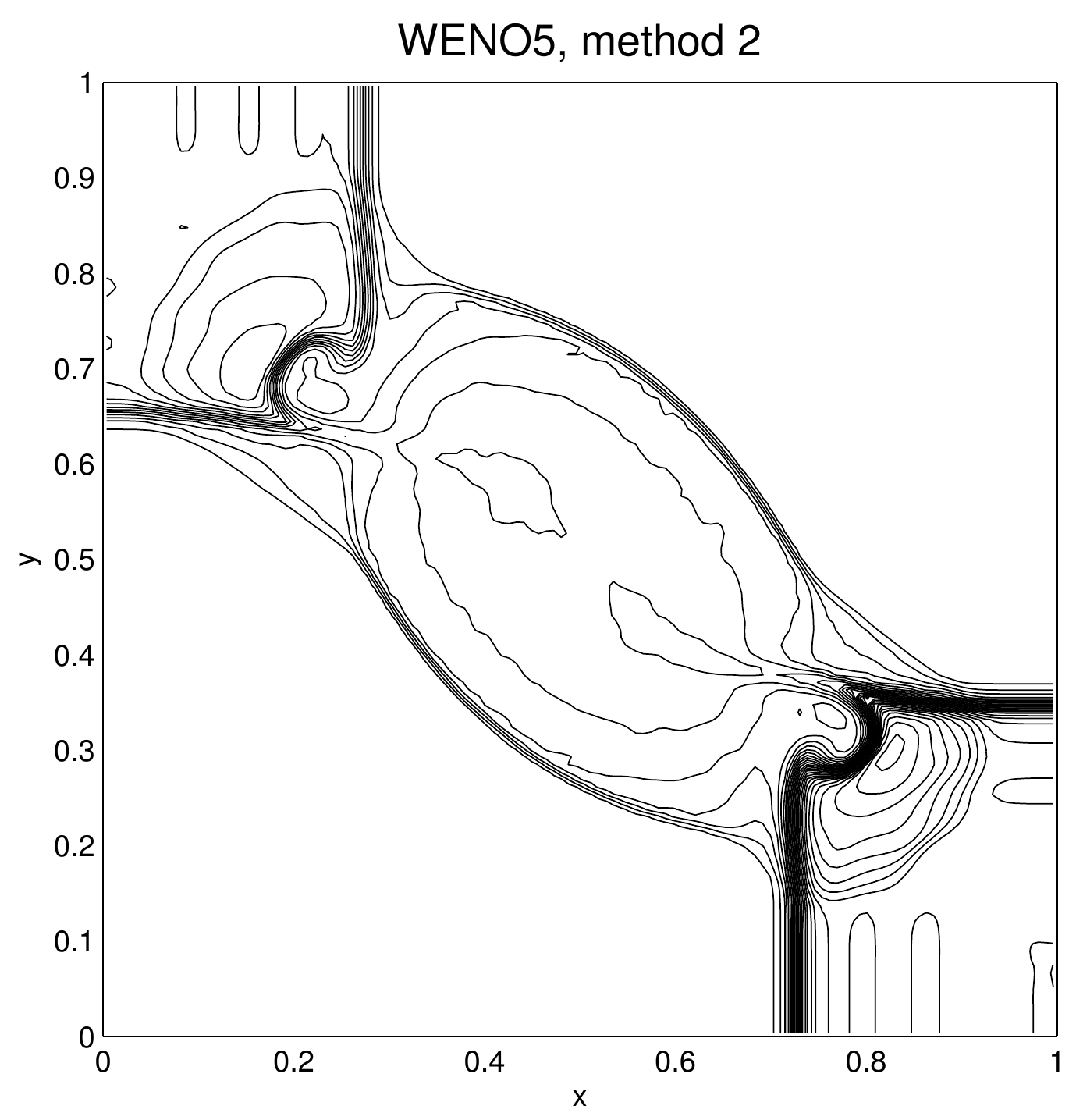}

\includegraphics[width=0.45\textwidth]{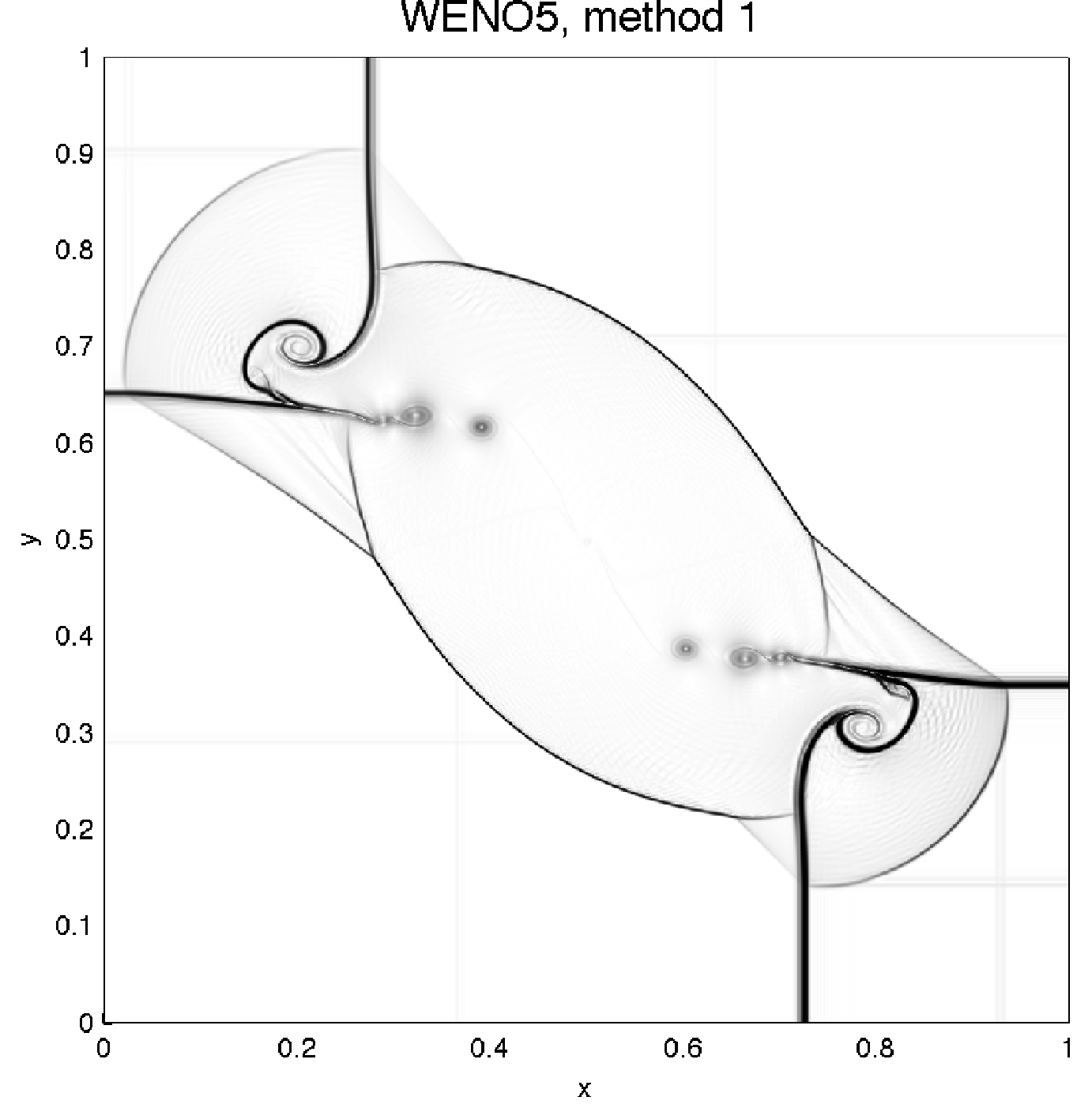}\hfil
\includegraphics[width=0.45\textwidth]{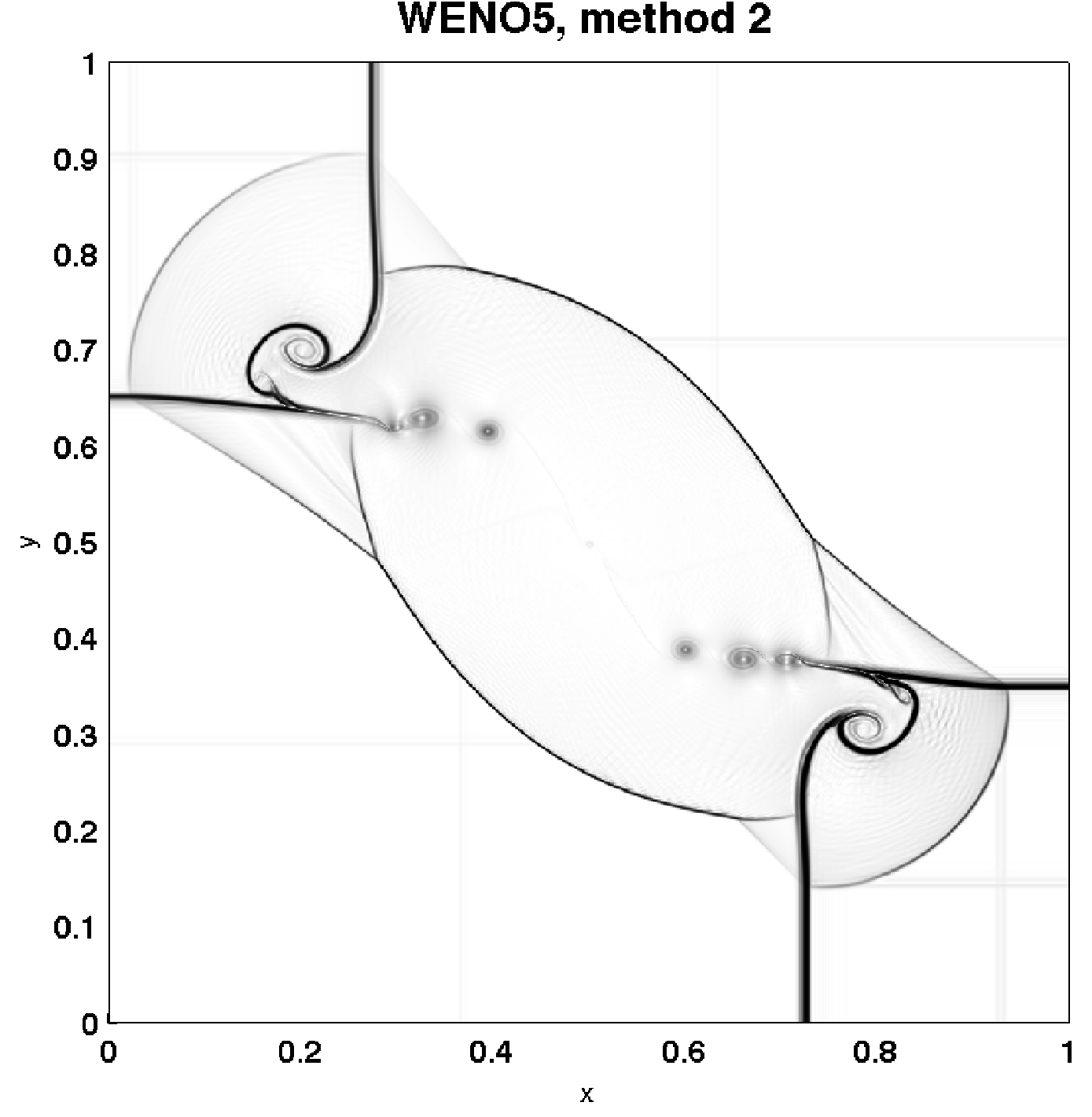}
\caption{\label{figure:2drp-1} Comparison of method 1 and method 2 for a standard
  two-dimensional Riemann problem on a grid with (top) $128 \times
  128$ and (bottom) $1024 \times 1024$ grid points.
Here we used WENO5 with RK5 and  the Roe Riemann solver. We show
contour lines of the computed solutions on coarse grids and schlieren
plots for the more resolved computations.}
\end{figure}

\begin{figure}[htb!]
\includegraphics[width=0.45\textwidth]{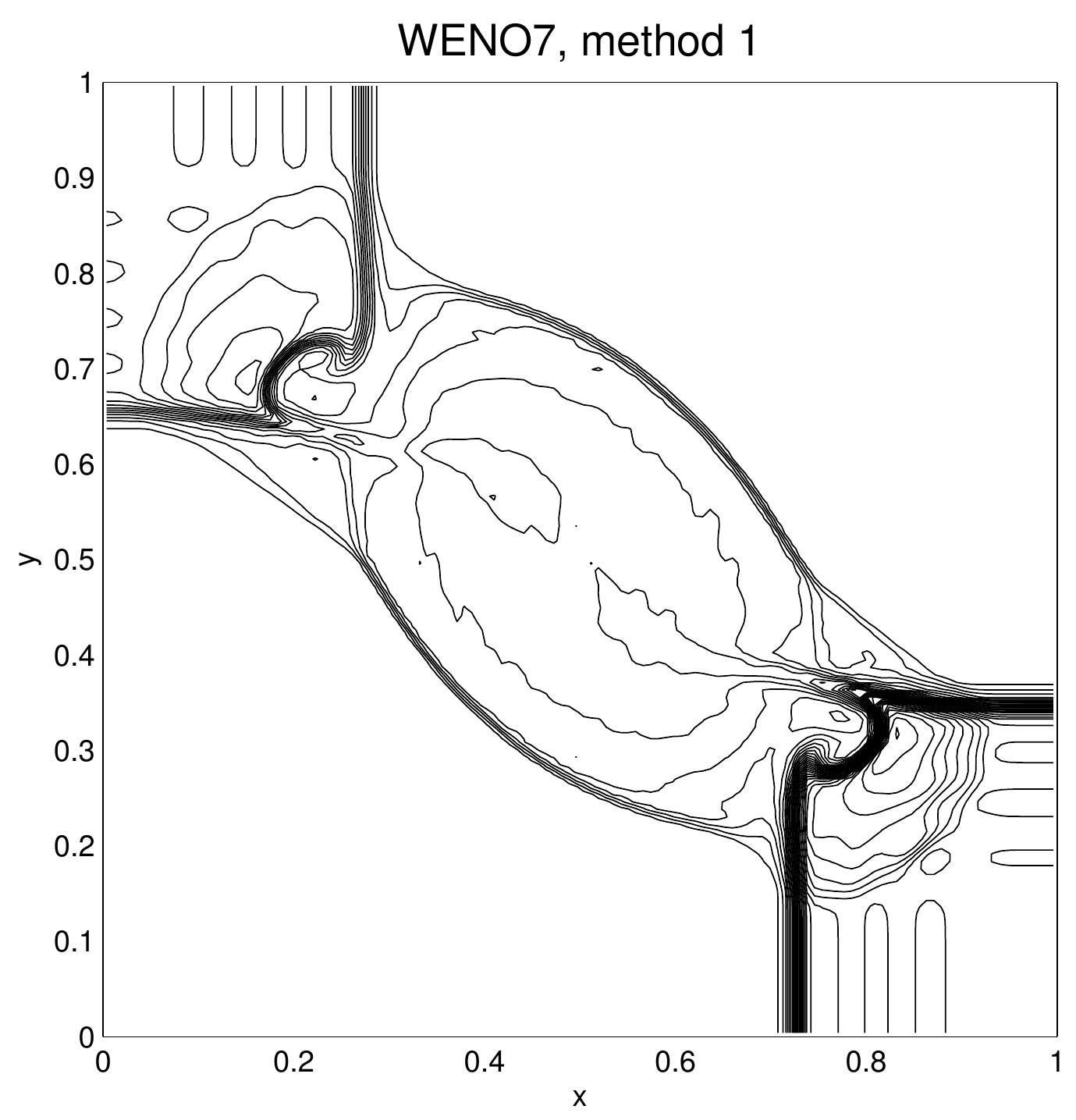}\hfil
\includegraphics[width=0.45\textwidth]{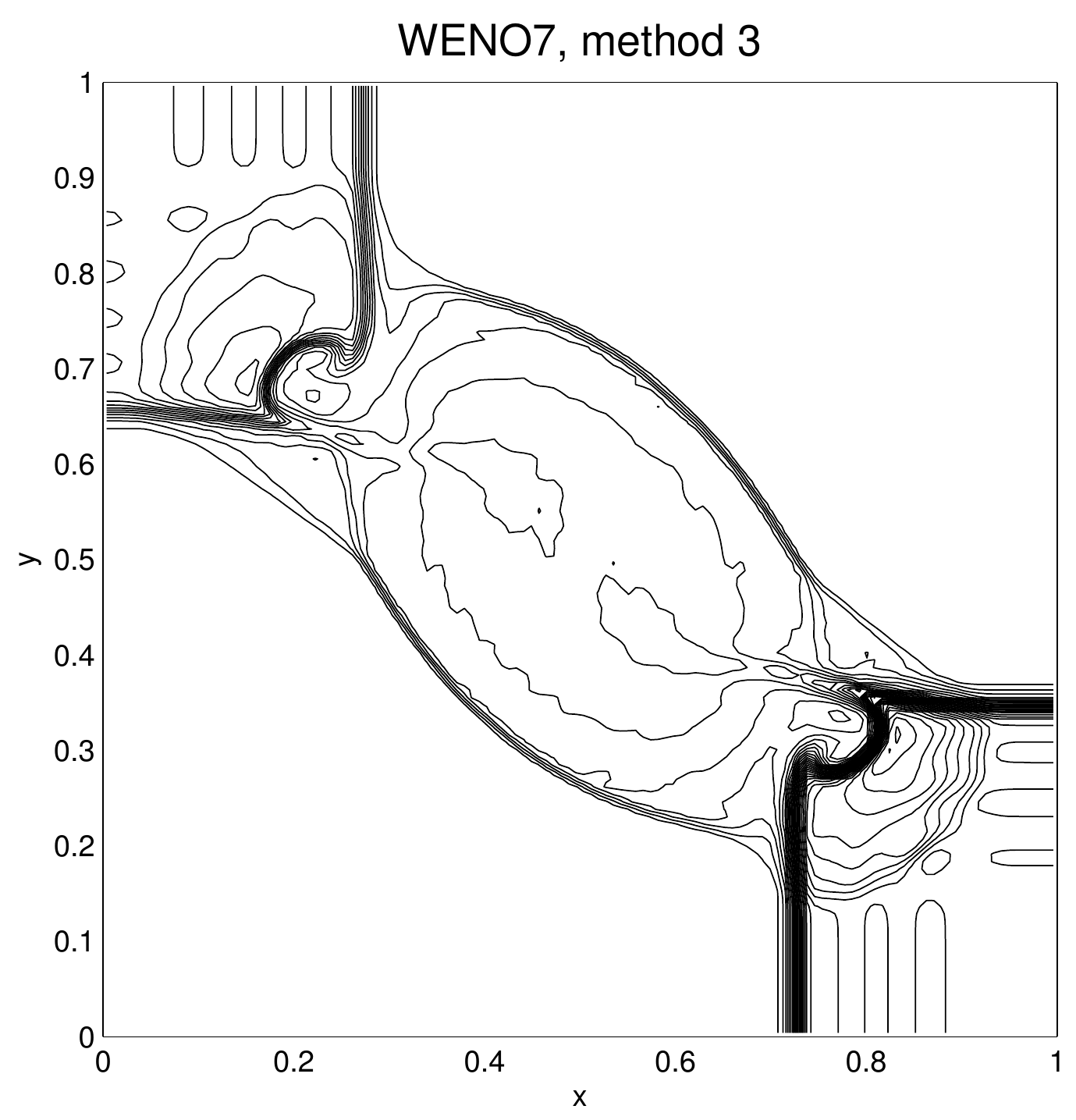}

\includegraphics[width=0.45\textwidth]{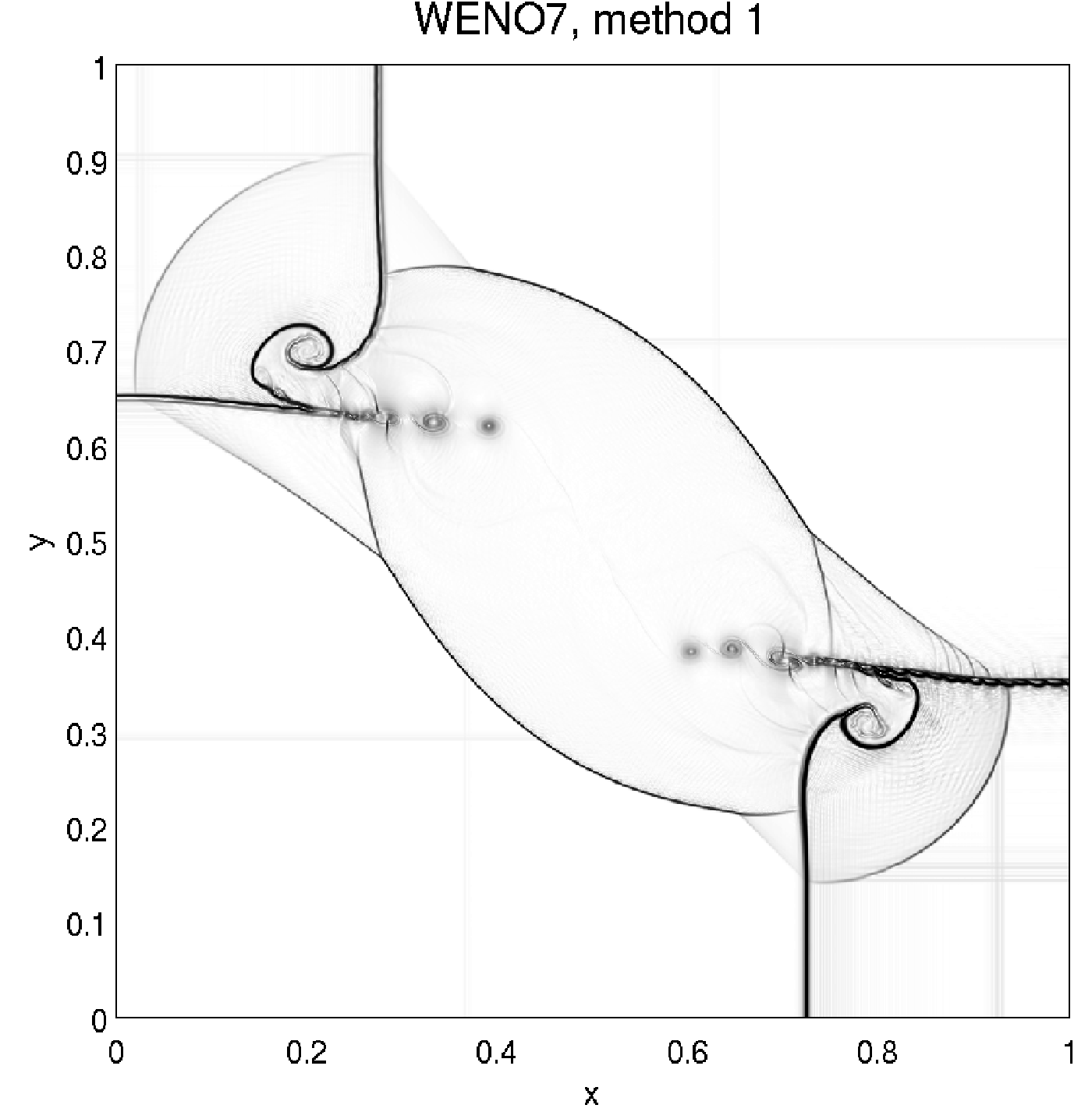}\hfil
\includegraphics[width=0.45\textwidth]{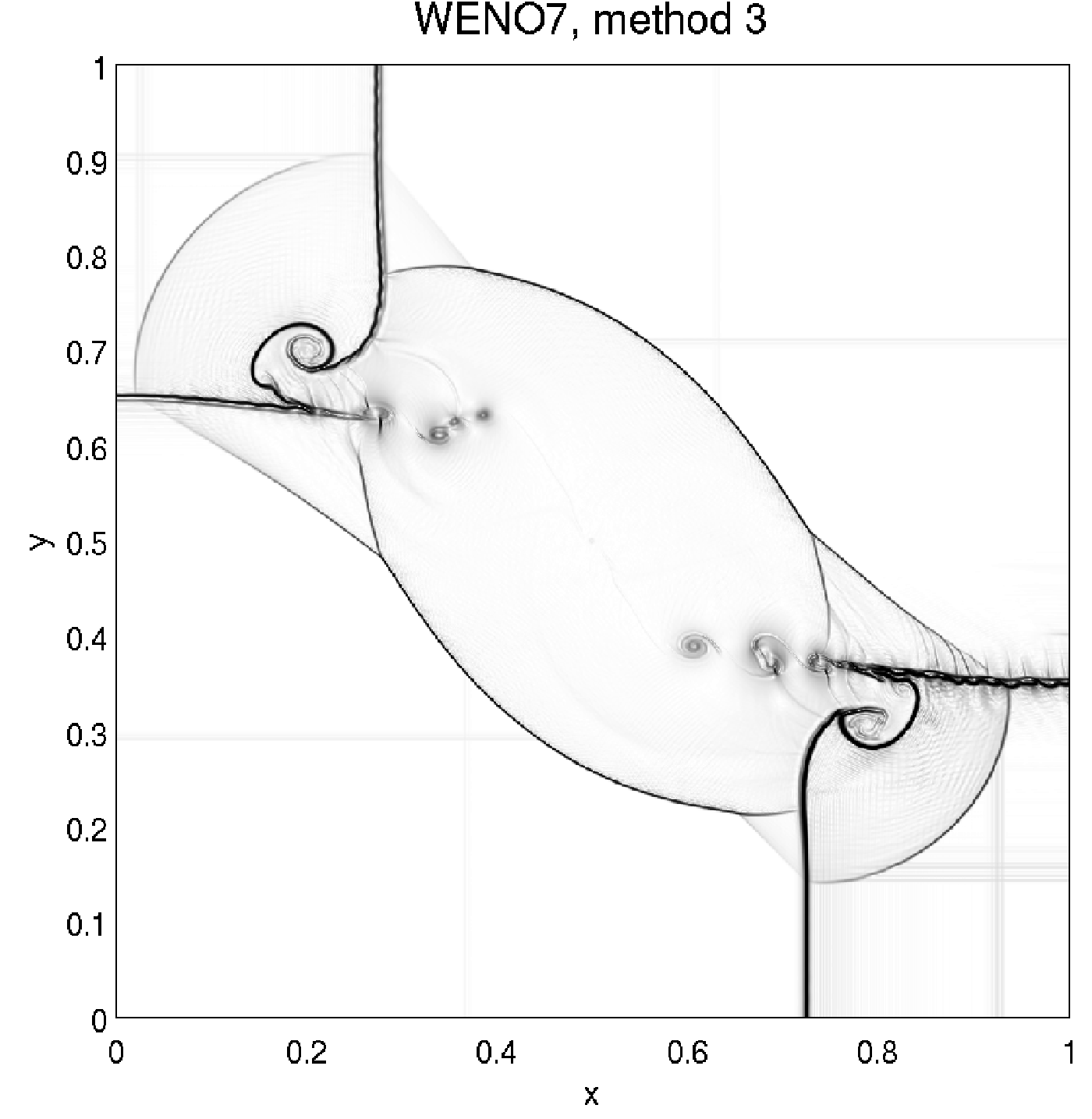}
\caption{\label{figure:2drp-2} Comparison of method 1 and method 3 for a standard
  two-dimensional Riemann problem on a grid with (top) $128 \times
  128$ and (bottom) $1024 \times 1024$ grid points.
Here we used WENO7 with RK7 and  the Roe Riemann solver. }
\end{figure}

In Figures \ref{figure:2drp-1} and \ref{figure:2drp-2} we show
results for \cite[configuration 5]{article:Schultz-Rinne1993}. 
The results obtained by the standard dimension--by--dimension approach
and our modified methods compare very well.
This is in agreements with observations reported in
\cite{article:ZZS2011}, where it was observed that for problems with
discontinuities  the classical dimension--by--dimension approach 
gives as good results
as the more expensive formally higher order accurate method. 
One could easily construct finite volume methods, which only use the
modified higher-order accurate update in regions where the solution is
smooth (and the modification is justified). 
Note that we used the WENO-JS instead of the WENO-Z
reconstruction (see Appendix \ref{appendix:s2}), since this
reconstruction produced 
slightly fewer oscillations behind the discontinuities.

\section{A high-order WENO finite volume method 
 for the equations of ideal magnetohydrodynamics}
Finally, we apply the modified WENO method to a more complex
application, namely the approximation of the 3d ideal
magnetohydrodynamic (MHD) equations.

The ideal MHD equations can be written in the form 
\begin{equation}\label{eqn:MHD}
\begin{split}
 \frac{\partial}{\partial t}
  \left( \begin{array}{c}
    \rho \\ \rho {\bf u} \\ {E} \\ {\bf B}
  \end{array} \right) +  \nabla \cdot \left(
  \begin{array}{c} \rho {\bf u} \\ \rho {\bf u}  {\bf u} + \left( {p} + \frac{1}{2}
  \| {\bf B} \|^2  \right) {\mathbb I}
    - {\bf B}  {\bf B} \\
    {\bf u} \left({E} + {p} + \frac{1}{2} \| {\bf B} \|^2 \right) - 
    	{\bf B} \left({\bf u} \cdot {\bf B} \right)
    \\ {\bf u}  {\bf B} - {\bf B}  {\bf u}
  \end{array} \right) & = 0, \\
  \nabla \cdot {\bf B} & = 0,
\end{split}
\end{equation}
where $\rho$, $\rho {\bf u}$ and ${E}$ are the total mass, momentum and
energy densities, and ${\bf B}$ is the magnetic field. The thermal pressure,
$p$, is related to the conserved quantities through the ideal gas law
\begin{equation}
\label{eqn:eos}
p = (\gamma -1) \left( {E} - \frac{1}{2} \| {\bf B} \|^2 - \frac{1}{2} \rho
  \| {\bf u} \|^2 \right),
\end{equation} 
where $\gamma = 5/3$ is the ideal gas constant. 

It is well known, that numerical methods for the multidimensional MHD
equations must control errors in the discrete divergence of the
magnetic field. One possibility to do this, is by using so-called
constrained transport (CT) methods. Here we use an approach which was
recently developed by Helzel, Rossmanith and
Taetz \cite{article:HRT2013,article:HRT2011}, and which is based on earlier work 
by Rossmanith \cite{article:R2006}.  
A FD-WENO method for the ideal MHD equations, using this same kind of
constrained transport, was recently proposed by
Christlieb, Rossmanith and Tang \cite{article:CRT2013}.

Since ${\bf B}$ is divergence free, we can set ${\bf B} = \nabla
\times {\bf A}$, where ${\bf A} \in \mathbb{R}^3$ is the magnetic potential.
Inserting this relation in the last line of the MHD equations, we
derive an evolution equation for the magnetic potential 
\begin{equation}\label{eqn:MHD-A}
 \partial_t {\bf A} + (\nabla \times {\bf A})
  \times {\bf u}  = - \nabla \phi.
\end{equation}
Here $\phi$ is an arbitrary scalar function. Different choices of
$\phi$ represent different {\em gauge condition} choices as explained
in  \cite{article:HRT2011}.  We use the so-called {\em Weyl gauge},
which means that we set $\nabla \phi = 0$ in (\ref{eqn:MHD-A}). The
resulting evolution equation for the magnetic potential can be written
in the form 
\begin{equation} \label{eqn:MHD-A-ql}
{\bf A}_t + N_1({\bf u}) {\bf A}_x + N_2({\bf u}) {\bf A}_y + N_3({\bf
  u}) {\bf A}_z = 0,
\end{equation}
with 
\begin{equation}\label{eqn:MHD-N}
N_1 = \left( \begin{array}{ccc}
0 & -u^2 & -u^3\\
0 & u^1 & 0\\
0 & 0 & u^1\end{array}\right), \ N_2 = \left( \begin{array}{ccc}
u^2 & 0 & 0\\
-u^1 & 0 & -u^3\\
0 & 0 & u^2\end{array}\right), \ N_3 = \left( \begin{array}{ccc}
u^3 & 0 & 0\\
0 & u^3 & 0\\
-u^1 & -u^2 & 0\end{array}\right).
\end{equation}
The system (\ref{eqn:MHD-A-ql}) with matrices of the form (\ref{eqn:MHD-N}) is weakly
hyperbolic, i.e.\ the matrix $N({\bf n}) = n^1 N_1({\bf u}) + n^2
N_2({\bf u}) + n^3 N_3({\bf u})$ has real eigenvalues for all ${\bf n}
\in S^2$, but there are directions for which $N({\bf n})$ fails to have a
complete set of right eigenvectors, see \cite{article:HRT2011}.   

To describe the general form of the constrained transport algorithm,
we introduce the notation
\begin{equation}
Q_{MHD}'(t) = {\cal L}_1 (Q_{MHD}(t)),
\end{equation}
for the semi-discrete form of the MHD equations. Here $Q_{MHD}(t)$
represents the grid function at time $t$ consisting of all
cell-averaged values of the conserved quantities from the MHD equation
(\ref{eqn:MHD}). Analogously, we introduce
\begin{equation}
Q_{\bf A}'(t) = {\cal L}_2 (Q_{\bf A} (t), Q_{MHD}(t)),
\end{equation}
to describe the semi-discrete form for the evolution equation of the
magnetic potential. Note that the evolution of the potential depends
on the velocity field, which we take to be as given function from the
solution step of the MHD equations.

To simplify notation, we present the numerical method using forward
Euler time-stepping. 
\begin{itemize}
\item[0.] Start with $Q_{MHD}^n$ and $Q_{\bf A}^n$ (i.e.\ the solution
  from the previous time step).
\item[1.] Update without regard of the divergence-free condition on
  the magnetic field, to obtain $Q_{MHD}^{*}$  and $Q_{\bf
    A}^{n+1}$:
\begin{eqnarray}
Q_{MHD}^* & = & Q_{MHD}^n + \Delta t {\cal L}_1
(Q_{MHD}^n) \label{eqn:MHD-1}\\
Q_{\bf A}^{n+1} & = & Q_{\bf A}^n + \Delta t {\cal L}_2 (Q_{\bf A}^n,
Q_{MHD}^n) \label{eqn:MHD-2}
\end{eqnarray}

\item[2.] Correct the magnetic field components $Q_{MHD}^{*}$ by
  the divergence-free values
${\bf B}^{n+1} = \nabla \times Q_{\bf A}^{n+1}$. Set $Q_{MHD}^{n+1} =
\left( \rho^{n+1}, \rho {\bf u}^{n+1}, E^{n+1}, {\bf B}^{n+1} \right)$. 
\end{itemize}

In Step 1, update (\ref{eqn:MHD-1}), we use a straight forward extension to the
three-dimensional case of our modified WENO method for hyperbolic
partial differential equations in
divergence form. Here we used the $5$th order WENO-Z method with a
correction that leads to fourth order accurate flux functions, i.e.\
method 2.     

In Step 1, update (\ref{eqn:MHD-2}), we used a three-dimensional
extension of our method from Section
\ref{section:s5}, to update the evolution equation for the magnetic
potential. Note that due to the weak hyperbolicity of (\ref{eqn:MHD-A-ql}), the
fluctuations ${\bf A}^\pm \Delta Q_{i+\frac{1}{2},j,k}$, ${\cal B}^\pm \Delta Q_{i,j+\frac{1}{2},k}$ and
${\cal C}^\pm \Delta Q_{i,j,k+\frac{1}{2}}$ can not be computed using an eigenvector
decomposition of the jump in $Q_{\bf A}$ at grid cell
interfaces. Instead, we computed the fluctuations using an approach
based on the idea of path conservative methods, as explained in \cite{article:HRT2013}.  

In Step 2, we compute ${\bf B}^{n+1}=(B^1,B^2,B^2)$ from the cell average values of
$Q_{\bf A}^{n+1} = (A^1,A^2,A^3)$, using the formulas
\begin{equation}
\begin{split}
B_{i,j,k}^1  =  &\frac{1}{12 \Delta y} \left(A_{i,j-2,k}^3 - 8
  A_{i,j-1,k}^3 + 8 A_{i,j+1,k}^3 - 
  A_{i,j+2,k}^3 \right) \\
& - \frac{1}{12 \Delta z} \left( A_{i,j,k-2}^2
 -8 A_{i,j,k-1}^2 + 8 A_{i,j,k+1}^2 - 
    A_{i,j,k+2}^2 \right)  \\
B_{i,j,k}^2 = & \frac{1}{12 \Delta z} \left( A_{i,j,k-2}^1 -
  8 A_{i,j,k-1}^1 + 8 A_{i,j,k+1}^1 - 
  A_{i,j,k+2}^1 \right) \\
& - \frac{1}{12 \Delta x} \left( A_{i-2,j,k}^3 - 8
  A_{i-1,j,k}^3 + 8 A_{i+1,j,k}^3 - 
  A_{i+2,j,k}^3 \right)\\
B_{i,j,k}^3 = & \frac{1}{12 \Delta x} \left( A_{i-2,j,k}^2 -
  8 A_{i-1,j,k}^2 + 8 A_{i+1,j,k}^2 - 
  A_{i+2,j,k}^2 \right) \\
& - \frac{1}{12 \Delta y} \left( A_{i,j-2,k}^1 -
    8 A_{i,j-1,k}^1 + 8 A_{i,j+1,k}^1 -
     A_{i,j+2,k}^1 \right). 
\end{split}
\end{equation}
This is a fourth order accurate approximation of  cell averaged values
of $\nabla \times {\bf A}$
using cell averaged values of the magnetic potential. 

We tested the new finite volume CT method for the 3d 
smooth Alfv\'en wave problem. The initial data and the computational
domain for this problem are described in \cite{article:HRT2011}. For
our computations we used WENO-Z5 + RK5. 
In Table \ref{table:MHD}, we show results of a numerical convergence study. 
There we compare the constrained transport method with our modified WENO method implemented in
method 2 with the simple dimension-by-dimension approach implemented
in method 1. As expected, we observe second order convergence rates for
method 1. Method 2 converges with fourth order.
\begin{table}[htb!]
  \centerline{
    \begin{tabular}{|c||c|c|c|c|}
      \hline
      & \multicolumn{2}{c|}{method 1} 
      & \multicolumn{2}{c|}{method 2} \\
     grid & $\|\rho-\rho_{ref}\|_1$ & EOC 
           & $\|\rho-\rho_{ref}\|_1$ & EOC  \\
\hline    
16x32x32 & 7.02906d-4 & & 7.81748d-4 & \\
32x64x64 & 6.23454d-5 & {\bf 3.49} & 5.28206d-5 & {\bf 3.89} \\
64x128x128 & 8.02708d-6 & {\bf 2.96} & 3.36290d-6 & {\bf 3.97} \\
128x256x256 & 1.57417d-6 & {\bf 2.35} & 2.11027d-7 & {\bf 3.99} \\
256x512x512 & 3.68469d-7 & {\bf 2.09} & 1.31970d-8 & {\bf 4.00}\\ \hline 
    \end{tabular}}
\caption{\label{table:MHD}Convergence study for smooth Alfv\'en wave problem using
  WENO-Z5+RK5 on a three-dimensional Cartesian grid.}
\end{table}

\section*{Conclusions}
We have presented a simple modification of the popular
dimension--by--dimen\-sion WENO method for Cartesian grids, which
retains the full order of accuracy of the corresponding
one--dimensional method. Our approach is based on a transformation of
interface values and point values of the conserved quantities and
numerical flux functions. 
For the popular WENO5 method, the simplest modification, which we call
method 2, gave very good results. For even higher order methods, such
as WENO7, we suggest to use the modification according to method 3.
 
Our method is an alternative to the 
multi--dimensional WENO finite volume method on Cartesian grids used
previously
\cite{article:SHS2002,article:TT2004,article:ZZS2011}. Compared to
theirs, our approach is less expensive. While the computing time for the method in  \cite{article:ZZS2011} is about
3.3-4.3 times that of the standard dimension--by--dimension approach,
our method requires only about 1.1-1.2  times the computing time of the
standard method.

For our considerations we always used the simplest version of the
finite volume WENO
method, where WENO reconstruction is performed component--wise for the
conserved quantities. Better results can often be obtained by applying
WENO reconstruction to primitive variables or characteristic
variables. The modifications suggested in this
paper  can also be introduced for such methods with only small modifications.

To simplify the notation, we have presented our improved versions of
the WENO method for two--dimensional problems. An extension to the
three-dimensional case is straight forward and the relevant formulas
for a transformation between cell averaged values and point values were
presented in Section \ref{sec:extens-high-dimens}. We have also applied 
three-dimensional versions of the proposed methods in the framework of
unstaggered constrained transport methods for the MHD equations.  

Together with J\"urgen Dreher, we are currently developing an AMR
version of our method.

\begin{appendix}
\section{Explicit Runge--Kutta methods}\label{app:RK4-6}
For the temporal discretization we use explicit Runge--Kutta methods of
order 5 and 7, respectively. After discretizing the PDE in space, we obtain a system
of ordinary differential equations of the general form
\begin{equation}
\frac{d}{dt} Q(t) = {\mathcal L}(Q(t)),
\end{equation}
where $Q(t)$ is a grid function of cell average values of the
conserved quantities at time $t$.
We discretize the resulting ode system using Runge--Kutta
methods of order five and seven. The methods are described by the Butcher
tableaus in Tables \ref{table:RK5}-\ref{table:RK7}.

\begin{table}[htb]
\begin{tabular}{c|c c c c c c}
$0$ & & & & & & \\[0.1cm]
$\frac{1}{4}$ & $\frac{1}{4}$ & & & & & \\[0.1cm]
$\frac{1}{4}$ & $\frac{1}{8}$ & $\frac{1}{8}$ & & & & \\[0.1cm]
$\frac{1}{2}$ & $0$ & $-\frac{1}{2}$ & $1$ & & &\\[0.1cm]
$\frac{3}{4}$ & $\frac{3}{16}$ & $0$ & $0$ & $\frac{9}{16}$ & &
\\[0.1cm]
$1$ & $-\frac{3}{7}$ & $\frac{2}{7}$ & $\frac{12}{7}$ &
$-\frac{12}{7}$ & $\frac{8}{7}$ & \\[0.1cm]
\hline

& $\frac{7}{90}$ & $0$& $\frac{32}{90}$ & $\frac{12}{90}$ &
$\frac{32}{90}$ & $\frac{7}{90}$
\end{tabular}
\caption{\label{table:RK5}Butcher tableau of the fifth order accurate
  Runge--Kutta method from \cite{article:RI2012}.}
\end{table}

\begin{figure}[htb!]
\includegraphics[width=0.45\textwidth]{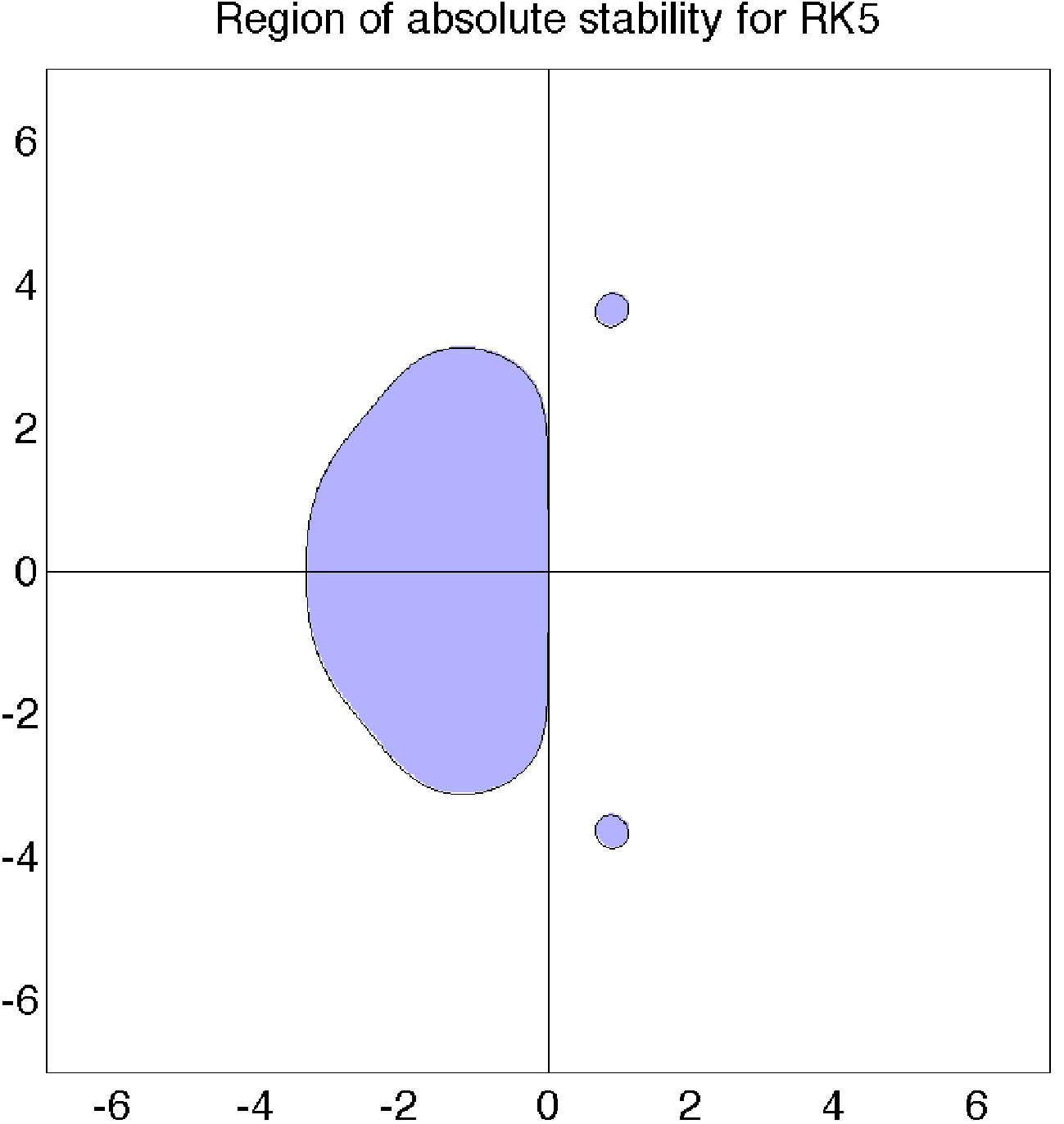}\hfill
\includegraphics[width=0.45\textwidth]{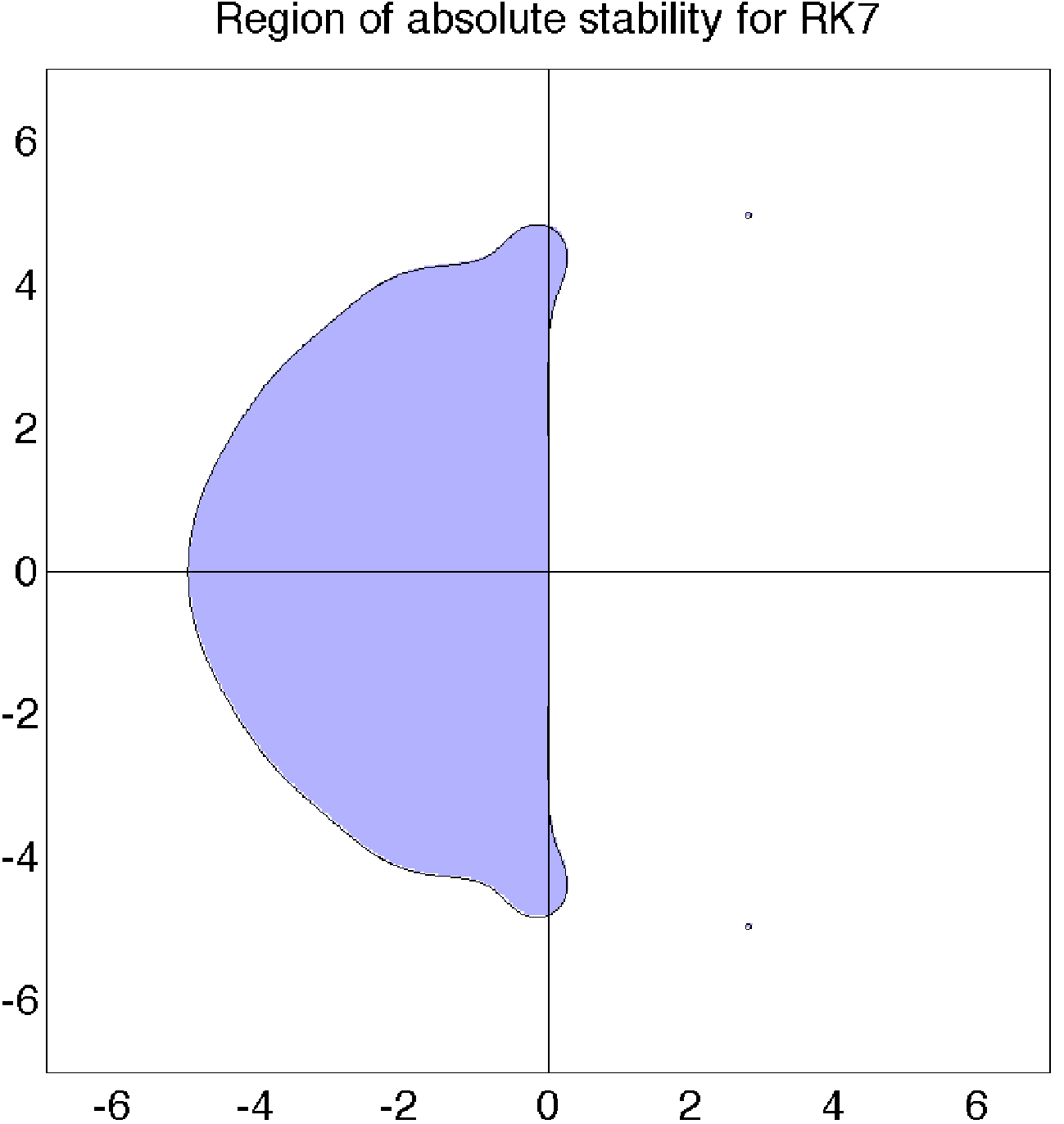}
\caption{\label{figure:rk-stab}
  Regions of absolute stability for RK5 and RK7.}   
\end{figure}

\begin{table}[htb]
\begin{tabular}{c|c c c c c c c c c c c}
$0$ & & & & & & & & & & &\\[0.1cm]
$\frac{2}{27}$ & $\frac{2}{27}$ & & & & & & & & & &\\[0.1cm]
$\frac{1}{9}$ & $\frac{1}{36}$ & $\frac{1}{12}$ & & & & & & & & &
\\[0.1cm]
$\frac{1}{6}$ & $\frac{1}{24}$ & $0$ & $\frac{1}{8}$ & & & & & &
& &\\[0.1cm]
$\frac{5}{12}$ & $\frac{5}{12}$ & $0$ & $-\frac{25}{16}$ &
$\frac{25}{16}$ & & & & & & &\\[0.1cm]
$\frac{1}{2}$ & $\frac{1}{20}$ & $0$ & $0$ & $\frac{1}{4}$
&$\frac{1}{5}$ & & & & & &\\[0.1cm]
$\frac{5}{6}$ & $-\frac{25}{108}$ & $0$ & $0$ & $\frac{125}{108}$ &
$-\frac{65}{27}$ & $\frac{125}{54}$ & & & & &\\[0.1cm]
$\frac{1}{6}$ & $\frac{31}{300}$ & $0$ & $0$ & $0$ & $\frac{61}{225}$
& $-\frac{2}{9}$ & $\frac{13}{900}$ & & & &\\[0.1cm]
$\frac{2}{3}$ & $2$ & $0$ & $0$ & $-\frac{53}{6}$ & $\frac{704}{45}$ &
$-\frac{107}{9}$ & $\frac{67}{90}$ & $3$ & & &\\[0.1cm]
$\frac{1}{3}$ & $-\frac{91}{108}$ & $0$ & $0$ & $\frac{23}{108}$ &
$-\frac{976}{135}$ & $\frac{311}{54}$ & $-\frac{19}{60}$ &
$\frac{17}{6}$ & $-\frac{1}{12}$ & &\\[0.1cm]
$1$ & $\frac{2383}{4100}$ & $0$ & $0$ & $-\frac{341}{164}$ &
$\frac{4496}{1025}$ & $-\frac{301}{82}$ & $\frac{2133}{4100}$ &
$\frac{45}{82}$ & $\frac{45}{164}$ & $\frac{18}{41}$ &\\[0.1cm]
\hline

& $\frac{41}{840}$ & $0$ & $0$ & $0$ & $0$ & $\frac{34}{105}$ &
$\frac{9}{35}$ & $\frac{9}{35}$ & $\frac{9}{280}$ & $\frac{9}{280}$ & $\frac{41}{840}$
\end{tabular}
\caption{\label{table:RK7}Butcher tableau of a seventh order accurate Runge--Kutta method
  from \cite{article:Fehlberg1969}.}
\end{table}
In Figure \ref{figure:rk-stab} we show the stability regions of the
two different Runge-Kutta methods used in this paper.

\section{Spatial reconstruction of interface values}\label{appendix:s2}
In this appendix we give the formulas for the spatial reconstruction
of interface averaged values of the conserved quantity used in our
implementation of the WENO method. We present the
formulas for the reconstruction in the $x$--direction.
This is based on a description of
WENO methods in \cite{article:BS2000,article:DB2013,article:Shu2009}. 
\subsection{$5$th order accurate WENO reconstruction}
At grid cell interfaces we compute averaged values of the conserved
quantities 
\begin{equation}\label{eqn:app2-1}
\begin{split}
Q_{i\pm\frac{1}{2},j}^{\mp} & = w_1^{\mp} Q_{i\pm\frac{1}{2},j}^{(1\mp)} + w_2^{\mp}
Q_{i\pm\frac{1}{2},j}^{(2\mp)} + w_3^{\mp} Q_{i\pm\frac{1}{2},j}^{(3\mp)},
\end{split}
\end{equation} 
with 
\begin{equation}
\begin{split}
Q_{i+\frac{1}{2},j}^{(1-)} & = \frac{1}{3} Q_{i-2,j} - \frac{7}{6}
Q_{i-1,j} + \frac{11}{6} Q_{i,j}, \quad Q_{i-\frac{1}{2},j}^{(1+)} =
-\frac{1}{6} Q_{i-2,j} + \frac{5}{6} Q_{i-1,j} + \frac{1}{3} Q_{i,j}\\
Q_{i+\frac{1}{2},j}^{(2-)}  & = -\frac{1}{6} Q_{i-1,j} + \frac{5}{6}
Q_{i,j} + \frac{1}{3} Q_{i+1,j}, \quad
Q_{i-\frac{1}{2},j}^{(2+)} = \frac{1}{3} Q_{i-1,j} + \frac{5}{6}
Q_{i,j} - \frac{1}{6} Q_{i+1,j} \\
Q_{i+\frac{1}{2},j}^{(3-)} & = \frac{1}{3} Q_{i,j} + \frac{5}{6}
Q_{i+1,j} - \frac{1}{6} Q_{i+2,j}, \quad 
Q_{i-\frac{1}{2},j}^{(3+)} = \frac{11}{6} Q_{i,j} - \frac{7}{6}
Q_{i+1,j} + \frac{1}{3} Q_{i+2,j}. \\
\end{split}
\end{equation}
The coefficients 
$w_1^\pm,\ldots, w_3^\pm$ 
in (\ref{eqn:app2-1}) depend
on the  local solution structure.
In the WENO-Z method suggested by Don and Borges \cite{article:DB2013},  they have the form
\begin{equation}\label{eqn:app2-2}
w_j^\pm = \frac{\tilde{w}_j^\pm}{\sum_{i=1}^3 \tilde{w}_i^\pm}, \quad \mbox{ with
} \tilde{w}_j^\pm = \gamma_j^\pm \left( 1 + \left(\frac{\tau_5}{\beta_j +
    \epsilon} \right)^p \right),
\end{equation}
with $p=2$ and $j=1,\ldots,3$.
$\gamma_1^- = \gamma_3^+ = \frac{1}{10}$, $\gamma_2^- = \gamma_2^+
= \frac{3}{5}$, $\gamma_3^- = \gamma_1^+ = \frac{3}{10}$,
$\beta_j$ as described in \cite[Equation (2.9)]{article:Shu2009},
$\tau_5 = |\beta_1-\beta_3|$ and 
$\epsilon=\Delta x^4$.
The WENO-Z methods are constructed to recover the optimal spatial
order of convergence. For other high order WENO methods this may
depend stronger on the choice of parameters such as $\epsilon$, see
for example \cite{article:HAP2005}.   

The WENO-JS method is obtained by replacing the computation of
$\tilde{w}_j^\pm$ by the formula
\begin{equation}\label{eqn:app2-3}
\tilde{w}_j^\pm = \frac{\gamma_j^\pm}{(\epsilon + \beta_j)^2},
\quad j=1,2,3.
\end{equation}
Here the same values are used for $\gamma_j^\pm$ and
$\beta_j$, but the parameter $\epsilon$ is replaced by $\epsilon=10^{-6}$.

\subsection{$7$th order accurate WENO reconstruction}
Our seventh order spatial spatial reconstruction uses values $Q^{(1
  \mp)}, \ldots, Q^{(4\mp)}$ from \cite{article:BS2000} and computes
$Q_{i\pm \frac{1}{2},j}^\mp$ analogously to (\ref{eqn:app2-1}), with
weights of the same form 
\begin{equation}\label{eqn:app2-4}
w_j^\pm = \frac{\tilde{w}_j^\pm}{\sum_{i=1}^4 \tilde{w}_i^\pm}, \quad \mbox{ with
} \tilde{w}_j^\pm = \gamma_j^\pm \left( 1 + \left(\frac{\tau_7}{\beta_j +
    \epsilon} \right)^p \right),
\end{equation}
The $\beta$--terms are set to be equal to $IS_0^4, \ldots,
IS_3^4$ as defined in \cite[page 415]{article:BS2000}. 
For the $7th$ order WENO-Z
method we set
$\tau_7 = |\beta_1 + 3 \beta_2 - 3 \beta_3-\beta_4|$, $p=2$ and
$\epsilon = \Delta x^5$, see \cite{article:DB2013}. 
For the $7th$ order WENO-JS method we use 
\begin{equation}\label{eqn:app2-5}
\tilde{w}_j^\pm = \frac{\gamma_j^\pm}{(\epsilon + \beta_j)^2},
\quad j=1,\ldots,4.
\end{equation}
with $\epsilon = 10^{-10}$.

\end{appendix}

% Non-BibTeX users please use


\begin{thebibliography}{10}
\bibitem{article:BS2000}
D.S.~Balsara and C.-W.~Shu.
\newblock Monotonicity preserving weighted essentially non-oscillatory schemes
  with increasingly high order of accuracy.
\newblock {\em J.\ Comput.\ Phys.}, 160:405--452, 2000.

\bibitem{article:CA1993}
J.~Casper and H.L.~Atkins.
\newblock A finite-volume high-order {ENO} scheme for two-dimensional
hyperbolic systems.
\newblock {\em J.\ Comput.\ Phys.}, 106:62--76, 1993.

\bibitem{article:CRT2013}
A.J.~Christlieb, J.A.~Rossmanith and Q.~Tang.
\newblock Finite difference weighted essentially non-oscillatory
schemes with constrained transport for ideal magnetohydrodynamics.
\newblock preprint, 2013. 

\bibitem{article:DB2013}
W.-S.~Don and R.~Borges.
\newblock Accuracy of the weighted essentially non-oscillatory
conservative finite difference schemes.
\newblock {\em J.\ Comput.\ Phys.}, 250: 347--372, 2013.


\bibitem{article:Fehlberg1969}
E.~Fehlberg.
\newblock Klassische {R}unge--{K}utta-{F}ormeln f\"unfter und siebenter
  {O}rdnung mit {S}chrittweiten-{K}ontrolle.
\newblock {\em Computing}, 4:93--106, 1969.

\bibitem{article:HRT2013}
C.~Helzel, J.A.~Rossmanith, and B.~Taetz.
\newblock A high order unstaggered constrained transport method for the ideal
  magnetohydrodynamic equations based on the method of lines.
\newblock {\em SIAM J.\ Sci.\ Comput.}, 35:A623--A651, 2013.

\bibitem{article:HRT2011}
C.~Helzel, J.A.~Rossmanith, and B.~Taetz.
\newblock An unstaggered constrained transport method for the 3d ideal
magnetohydrodynamic equations.
\newblock {\em J.\ Comput.\ Phys.}, 230: 3803--3829, 2011.

\bibitem{article:HAP2005}
A.K.~Henrick, T.D.~Aslam, and J.M.~Powers.
\newblock Mapped weighted essentially non-oscillatory schemes:
Achieving optimal order near critical points.
\newblock {\em J.\ Comput.\ Phys.}, 207:542--567, 2005.

\bibitem{article:HS1999}
C.~Hu and C.-W.~Shu.
\newblock Weighted essentially non-oscillatory schemes on triangular meshes.
\newblock {\em J.\ Comput.\ Phys.}, 150:97--127, 1999.

\bibitem{article:KPL2013}
D.I.~Ketcheson, M.~Parsani, and R.J.~LeVeque.
\newblock High-order wave propagation algorithms for hyperbolic systems.
\newblock {\em SIAM J.\ Sci.\ Comput.}, 35:A351--A377, 2013.

\bibitem{book:RJL2002}
R.J.~LeVeque.
\newblock {\em Finite Volume Methods for Hyperbolic Problems}.
\newblock Cambridge University Press, 2002.

\bibitem{article:MC2011}
P.~McCorquodale and P.~Colella.
\newblock A high--order finite volume method for conservation laws on logically
  refined grids.
\newblock {\em Commun.\ App.\ Math. and Comp.\ Sci.}, 6:1--25, 2011.

\bibitem{article:Merriman2003}
B.~Merriman.
\newblock Understanding the Shu--Osher conservative finite difference
form.
\newblock {\em J.\ Sci.\ Comput.}, 19:309--322, 2003.

\bibitem{article:RI2012}
F.~Rabiei and F.~Ismail.
\newblock Fifth-order improved {R}unge-{K}utta methods with reduced
number of function evaluations.
\newblock {\em Australien Journal of Basic and Applied Sciences},
6:97--105, 2012.

\bibitem{article:R2006}
J.A.~Rossmanith.
\newblock An unstaggered, high-resolution constrained transport method
for magnetohydrodynamic flows.
\newblock {\em SIAM J.\ Sci.\ Comput.}, 28:1766--1797, 2006 

\bibitem{article:TTD2011}
P.~Tsoutsanis, V.A.~Titarev, and D.~Drikakis.
\newblock Weno schemes on arbitrary mixed--element unstructured meshes in three
  space dimension.
\newblock {\em J.\ Comput.\ Phys.}, 230:1585--1601, 2011.

\bibitem{article:Schultz-Rinne1993}
C.W.~Schultz-Rinne.
\newblock Classification of the Riemann problem for two dimensional gas
  dynamics.
\newblock {\em SIAM J.\ Math.\ Anal.}, 24:76--88, 1993.

\bibitem{article:SQC2011}
C.~Shen, J.M.~Qiu and A.~Christlieb.
\newblock Adaptive mesh refinement based on high order finite
difference WENO scheme for multi-scale simulations.
\newblock {\em J.\ Comput.\ Phys.}, 230:3780--3802, 2011.

\bibitem{article:SZ2010}
Y.~Shen and G.~Zha.
\newblock Improved seventh-order {WENO} schemes.
\newblock {\em AIAA Paper}, 2010-1451.

\bibitem{article:Shu2009}
C.-W.~Shu.
\newblock High order weighted essentially nonoscillatory schemes for convection
  dominated problems.
\newblock {\em SIAM Review}, 51:82--126, 2009.

\bibitem{article:SO88}
C.-W.~Shu and S.~Osher.
\newblock Efficient implementation of essentially non-oscillatory
shock-capturing schemes.
\newblock {\em J.\ Comput.\ Phys.}, 77:439--471, 1988.

\bibitem{article:SO89}
C.-W.~Shu and S.~Osher.
\newblock Efficient implementation of essentially non-oscillatory
shock-capturing schemes, II.
\newblock {\em J.\ Comput.\ Phys.}, 83:32--78, 1989.

\bibitem{article:SHS2002}
\newblock J.~Shi, C.~Hu and C.-W.~Shu.
\newblock A technique for treating negative weights in WENO schemes.
\newblock {\em J.\ Comput.\ Phys.}, 175:108--127,  2002.

\bibitem{article:TT2004}
\newblock V.A.~Titarev and E.F.~Toro.
\newblock Finite-volume WENO schemes for three-dimensional conservation laws.
\newblock {\em J.\ Comput.\ Phys.}, 201:238--260, 2004.

\bibitem{article:ZZS2011}
R.~Zhang, M.~Zhang, and C.-W.~Shu.
\newblock On the order of accuracy and numerical performance of two classes of
  finite volume {WENO} schemes.
\newblock {\em Commun.\ Comput.\ Phys.}, 9:807--827, 2011.



\end{thebibliography}
\end{document}